\documentclass[a4paper,twosided,12pt]{amsart}

\usepackage[latin1]{inputenc}
\usepackage{amsxtra}
\usepackage{amssymb}
\usepackage{amsmath,amsthm}
\usepackage[all]{xypic}
\usepackage{epsfig}
\xyoption{dvips}

\theoremstyle{definition}
\newtheorem{ntn}{Notation}[subsection]
\newtheorem{dfn}[ntn]{Definition}
\newtheorem{rem}[ntn]{Remark}
\newtheorem{exa}[ntn]{Example}
\theoremstyle{plain}
\newtheorem*{que}{Question}
\newtheorem*{thmintr}{Theorem}

\newtheorem{lem}[ntn]{Lemma}

\newtheorem{prp}[ntn]{Proposition}
\newtheorem{thm}[ntn]{Theorem}
\newtheorem{cor}[ntn]{Corollary}
\newtheorem{claim}[ntn]{Claim}
\newtheorem{assmt}[ntn]{Assumption}
\theoremstyle{remark}

\DeclareMathAlphabet{\mathds}{U}{dsrom}{m}{n}
\DeclareMathAlphabet{\mathsc}{U}{rsfs}{m}{n}

\DeclareMathOperator{\id}{id}
\DeclareMathOperator{\Sch}{Sch}
\DeclareMathOperator{\Gal}{Gal}
\DeclareMathOperator{\Spec}{Spec}
\DeclareMathOperator{\Spf}{Spf}
\DeclareMathOperator{\Hg}{Hg}
\DeclareMathOperator{\MT}{MT}
\DeclareMathOperator{\Hm}{Hom}

\DeclareMathOperator{\Pic}{Pic}
\DeclareMathOperator{\End}{End}
\DeclareMathOperator{\NS}{NS}
\DeclareMathOperator{\DefF}{Def}
\DeclareMathOperator{\CSpin}{CSpin}
\DeclareMathOperator{\CSp}{CSp}

\DeclareMathOperator{\SO}{SO}
\DeclareMathOperator{\GL}{GL}

\DeclareMathOperator{\codim}{codim}

\newcommand{\Q}{\mathbb{Q}}
\newcommand{\N}{\mathbb{N}}
\newcommand{\Z}{\mathbb{Z}}
\newcommand{\F}{\mathbb{F}}

\newcommand{\R}{\mathbb{R}}
\newcommand{\A}{\mathbb{A}}
\newcommand{\C}{\mathbb{C}}
\newcommand{\Gm}{\mathbb{G}}
\newcommand{\s}{\mathbb{S}}
\newcommand{\Kg}{\mathbb{K}}

\renewcommand{\O}{\mathcal O}
\renewcommand{\L}{\mathcal L}

\newcommand{\M}{\mathcal M}

\newcommand{\X}{\mathcal X}

\renewcommand{\a}{\alpha}
\renewcommand{\l}{\lambda}
\renewcommand{\b}{\beta}
\newcommand{\G}{\Gamma}

\newcommand{\e}{\epsilon}

\newcommand{\Fk}{\mathcal F}
\newcommand{\Av}{\mathcal A}

\newcommand{\lr}{\rightarrow}

\title{Kuga-Satake Abelian Varieties in Positive Characteristic}
\author{Jordan Rizov}
\address{}
\email{danence@gmail.com}

\keywords{K3 surfaces, Moduli spaces, Kuga-Satake abelian varieties}
\subjclass[2000]{14J28, 14D22}

\begin{document}
\maketitle
\begin{abstract}
Kuga and Satake associate with every polarized complex K3 surface $(X,\L)$ a complex abelian variety called the Kuga-Satake abelian variety of $(X,\L)$. We use this construction to define morphisms between moduli spaces of polarized K3 surface with certain level structures and moduli spaces of polarized abelian varieties with level structure over $\C$. In this note we study these morphisms. We prove first that they are defined over finite extensions of $\Q$. Then we show that they extend in positive characteristic. In this way we give an indirect construction of Kuga-Satake abelian varieties over an arbitrary base. We also give some applications of this construction to canonical lifts of ordinary K3 surfaces.
\end{abstract}
\section*{Introduction}

When studying algebraic curves one constructs an abelian variety, called the Jacobian of the curve. The geometry of this abelian variety describes properties of the curve. Here, we consider a similar construction for K3 surfaces. Namely, we assign to every polarized K3 surface an abelian variety with certain properties, called its Kuga-Satake abelian variety. 

Let us explain briefly the construction of these varieties over $\C$ due to Kuga and Satake. Starting with a polarized complex K3 surface $(X,\L)$ one considers the second primitive Betti cohomology group
\begin{displaymath}
 P^2_B(X,\Z(1)) := c_1(\L)^\perp \subset H^2_B(X,\Z(1)).
\end{displaymath}
The orthogonal complement is taken with respect to the Poincar\'e pairing on $H^2_B(X,\Z(1))$. Using the polarized $\Z$-Hodge structure on $P^2_B(X,\Z(1))$ one defines a polarized $\Z$-Hodge structure of type $\{(1,0),(0,1)\}$ on the even Clifford algebra $C^+\bigl(P^2_B(X,\Z(1))\bigr)$. One might think of this construction as ``taking a square root of a Hodge structure''. Such a Hodge structure corresponds to a complex abelian variety $A$, called the Kuga-Satake abelian variety associated to $(X,\L)$. Using Kuga-Satake varieties one can deduce some properties of K3 surfaces, mostly of motivic nature, from the corresponding properties of abelian varieties.

At this point one may ask whether one can use this construction to define Kuga-Satake abelian varieties over subfields of $\C$. Or whether one can construct Kuga-Satake abelian schemes starting with families of polarized K3 surfaces. One can find some answers in \cite{D-K3} and \cite{A-HV}. One can go even further and ask whether one can define Kuga-Satake abelian varieties in positive characteristic. We combine all these questions in one, which was originally the motivation for our work.
\begin{que}
Can one define Kuga-Satake abelian varieties using only methods of algebraic geometry without making use of complex analytic constructions?
\end{que}
Up to isogeny a positive answer to this question can be found in Theorem 7.1 in \cite{A-MHC} and Theorem 1.5.1 in \cite{A-HV}. We refer also to Chapter 9 and 10.2.4 in \cite{Andre-Mot}. Starting with a polarized K3 surface $(X,\L)$, Y. Andr\'e constructs a ``motive'' which is isomorphic to the motive of the Kuga-Satake abelian variety of $(X,\L)$. The construction is purely algebro-geometric.

Here we solve a modification of this problem. Namely, we will be interested in a Kuga-Satake construction over an arbitrary base without putting any restriction on the ``methods''. The reason is that we use the existing transcendental construction as a starting point. We explain this in more detail.

P. Deligne gives an interpretation of the Kuga-Satake construction in terms of the adjoint representation homomorphism $\CSpin(2,19) \lr \SO(2,19)$ and the spin representation homomorphism $\CSpin(2,19) \hookrightarrow \CSp_{2g}$, where $g = 2^{19}$ (see \S\S 3 and 4 in~\cite{D-K3}). We consider the morphisms, between the Shimura varieties associated to the groups $\CSpin(2,19),\ \SO(2,19)$ and $\CSp_{2g}$, defined by the adjoint and the spin representations. Putting together these maps and the results of \cite{JR-CMK3}, for every $n \geq 3$, we define a Kuga-Satake morphism
\begin{displaymath}
f^{ks}_{d,a,n,\gamma,E_n} \colon \Fk_{2d,n^{\rm sp}}\otimes E_n \lr \Av_{g,d',n}\otimes E_n
\end{displaymath}
where $\Fk_{2d,n^{\rm sp}}$ is the moduli space of K3 surfaces with a primitive polarization of degree $2d$ and a spin level $n$-structure (see Section 6 in \cite{Riz-MK3}), and $E_n$ is a finite abelian extension of $\Q$. The morphism $f^{ks}_{d,a,n,\gamma,E_n}$ assigns to every primitively polarized complex K3 surface with a spin level $n$-structure its associated Kuga-Satake abelian variety plus extra data (a polarization and a level $n$-structure). In this way, the first step of our program is completed. Next we show that $f^{ks}_{d,a,n,\gamma,E_n}$ extends over an open part of $\Spec(\O_{E_n})$ where $\O_{E_n}$ is the ring of integers in $E_n$. More precisely, we prove the following statement.
\begin{thmintr}
Let $d,n \in \N$ and suppose that $n\geq 3$. Then the Kuga-Satake morphism $f^{ks}_{d,n,a,\gamma,E_n} \colon \Fk_{2d,n^{\rm sp},E_n} \lr \Av_{g,d',n,E_n}$ extends uniquely to a morphism
\begin{displaymath}
 f^{ks}_{d,a,n,\gamma} \colon \Fk_{2d,n^{\rm sp}}\otimes \O_{E_n}[1/N] \lr \Av_{g,d',n} \otimes \O_{E_n}[1/N]
\end{displaymath}
where $N = 2dd'nl$ and $l$ is the product of the prime numbers $p$ whose ramification index $e_p$ in $E_n$ is $\geq p-1$.
\end{thmintr}
The proof of this theorem is based on a result of G. Faltings on extension of abelian schemes.

We conclude with some applications. We show that the \'etale cohomology relations from \S 6.6.1 in \cite{D-K3} hold for the Kuga-Satake abelian varieties we construct. Then we focus our attention on the ordinary locus of $\Fk_{2d,n^{\rm sp}}\otimes \F_p$, where $p$ is a prime not dividing $N$. Suppose that $k$ is a finite field of characteristic $p$. One can easily see that $f^{ks}_{d,a,n,\gamma}$ maps an ordinary point $x = (X,\L,\nu)$ in $\Fk_{2d,n^{\rm sp}}\otimes \F_p (k)$ to an ordinary point $y = (A,\mu,\e)$ in $\Av_{g,d',n} \otimes \F_p (k)$. Denote by $x^{\rm can} = (X^{\rm can}, \L, \nu)$ the canonical lift of $X$ over $W(k)$ and by $y^{\rm can} = (A^{\rm can},\mu,\e)$ the canonical lift of $A$. We prove that $f^{ks}_{d,a,n,\gamma}(x^{\rm can}) = y^{\rm can}$. A straightforward corollary of this is that the restriction of the Kuga-Satake morphism to the ordinary locus of $\Fk_{2d,n^{\rm sp}}\otimes \F_p$ is quasi-finite.
\newline
\newline
{\bf Notations and conventions}
\newline
\newline
{\bf General.} We write $\hat \Z$ for the profinite completion of $\Z$. We denote by $\A$ the ring of ad\`eles of $\Q$ and by $\A_f = \hat \Z \otimes \Q$ the ring of finite ad\`eles of $\Q$. Similarly, for a number field $E$ we denote by $\A_E$ and $\A_{E,f}$ the ring of ad\`eles and the ring of finite ad\`eles of $E$.

If $A$ is a ring, $A \lr B$ a ring homomorphism then for any $A$-module ($A$-algebra etc.) $V$ we will denote by $V_B$ the $B$-module ($B$-algebra etc.) $V\otimes_A B$.

For a variety $X$ over $\C$ we will denote by $X^{\rm an}$ the associated analytic variety. For an algebraic stack $\mathcal F$ over a scheme $S$ and a morphism of schemes $S' \lr S$ we will denote by $\mathcal F_{S'}$ the product $\mathcal F\times_S S'$ and consider it as an algebraic stack over $S'$.

We use the notations established in \cite{Riz-MK3}. In particular, for a natural number $d$ we write $\Fk_{2d}$ for the Deligne-Mumford stack of K3 spaces with a primitive polarization of degree $2d$. It is a smooth stack over $\Spec(\Z[1/2d])$. See Theorem 4.7 in \cite{Riz-MK3} and \S 1.4.3 in \cite{JR-Thesis}. For $n \in \N$, $n \geq 3$, and a subgroup $\Kg$ of finite index in $\Kg_n$ we denote by $\Fk_{2d,\Kg}$ the smooth algebraic space over $\Spec(\Z[1/N])$ of K3 surfaces with a primitive polarization of degree $2d$ and a level $\Kg$-structure. If $\Kg$ is admissible, then we denote by $\Fk^{\rm full}_{2d,\Kg}$ the moduli space of K3 surfaces with a primitive polarization of degree $2d$ an a full level $\Kg$-structure. For details we refer to Section 6 in \cite{Riz-MK3} and Section 1.5 in \cite{JR-Thesis}.

Let $U$ be the hyperbolic plane and let $E_8$ be the \emph{positive} quadratic lattice associated to the Dynkin diagram of type $E_8$. Denote by $(L_0,\psi)$ the quadratic lattice $U^{\oplus 3} \oplus E_8^{\oplus 2}$. Further, let $(V_0,\psi_0)$ be the quadratic space $(L_0,\psi) \otimes_\Z \Q$. Further, we use the notations established in Section 2.1 of \cite{Riz-MK3}. 
\newline
\newline 
{\bf Algebraic groups.} A superscript ${}^0$ usually indicates a connected component for the Zariski topology. For an algebraic group $G$ will denote by $G^0$ the connected component of the identity. We will use the superscript ${}^+$ to denote connected components for other topologies.

For a reductive group $G$ over $\Q$ we denote by $G^{\rm ad}$ the adjoint group of $G$, by $G^{\rm der}$ the derived group of $G$ and by $G^{\rm ab}$ the maximal abelian quotient of $G$. We let $G(\R)_+$ denote the group of elements of $G(\R)$ whose image in $G^{\rm ad}(\R)$ lies in its identity component $G^{\rm ad}(\R)^+$, and we let $G(\Q)_+ = G(\Q) \cap G(\R)_+$. 

Let $V$ be a vector space over $\Q$ and let $G \hookrightarrow \GL(V)$ be an algebraic group over $\Q$. Suppose given a full lattice $L$ in $V$ (i.e., $L\otimes \Q = V$). Then $G(\Z)$ and $G(\hat \Z)$ will denote the abstract groups consisting of the elements in $G(\Q)$ and $G(\A_f)$ preserving the lattices $L$ and $L_{\hat \Z}$ respectively.
\newline
\newline
{\bf Acknowledgments}
\newline
\newline
This note contains the results of Chapter 4 of my Ph.D. thesis \cite{JR-Thesis}. I thank my advisors, Ben Moonen and Frans Oort for their help, their support and for everything I have learned from them. I would like to thank Ben Moonen for pointing out some mistakes in the earlier versions of the text and for his valuable suggestions. I thank the Dutch Organization for Research N.W.O. for the financial support with which my thesis was done.
\section{Extension of polarizations of abelian schemes}
In this section, we give some results on extension of polarizations of abelian schemes. We will use them to extend Kuga-Satake morphisms in positive characteristic. We fix a discrete valuation ring $R$ with field of fractions $K$ and residue field $k$.
\begin{lem}\label{pol}
Let $A$ be an abelian scheme over a discrete valuation ring $R$ and let $\l_K$ be a polarization of the generic fiber $A_K$ of $A$. Then $\l_K$ extends uniquely to a polarization of $A$.
\end{lem}
\begin{proof}
By Proposition 2.7 in Ch. 1 in \cite{CF-AV} $\l_K$ extends uniquely to a homomorphism $\l \colon A \lr A^t$ over $R$. It suffices to show that $2\cdot \l$ is a polarization. But $2\cdot \l = \varphi_{\mathcal M}$ where $\mathcal M = (\id_A,\l)^*\mathcal P_A$ and $\mathcal P_A$ is the Poincar\'e bundle on $A\times A^t$. We conclude by Corollary VIII 7 in \cite{R-AShGrSch} that $\M$ is relatively ample, hence $2\cdot \l$ is a polarization.
\end{proof}
\begin{lem}\label{polext}
Suppose given a locally noetherian, regular scheme $U$ and a dense open subscheme $V \subset U$ such that the codimension of $U \setminus V$ in $U$ is at least $2$. Let $A \lr U$ be an abelian scheme and let $\l_V$ be a polarization of $A_V \lr V$. Then $\l_V$ extends uniquely to a polarization $\l$ of $A \lr U$.
\end{lem}
\begin{proof}
Applying Proposition 2.7 of Ch. 1 in \cite{CF-AV} as in the proof of the previous lemma we see that $\l_V$ extends uniquely to an isogeny $\l \colon A \lr A^t$ over $U$.

By assumption there is an \'etale covering $\pi_V \colon \tilde V \lr V$ such that the pull-back $\l_{\tilde V} \colon A_{\tilde V} \lr A^t_{\tilde V}$ of $\l_V$ is equal to $\varphi_{\L_{\tilde V}}$ for an ample line bundle $\L_{\tilde V}$ on $A_{\tilde V}$. By the Zariski-Nagata purity theorem (see Cor. 3.3 of Exp. X in \cite{SGA1}), the morphism $\pi_V$ extends to an \'etale covering $\pi \colon \tilde U \lr U$. Let $j \colon \tilde V \lr \tilde U$ be the inclusion. Then the sheaf $\L := j_*\L_{\tilde V}$ is a line bundle (cf. Lemma 6.2 in Ch. V in \cite{CF-AV}). The isogenies $\l_{\tilde V}$ and $\varphi_{\L_{\tilde V}}$ coincide so by the unicity part of Proposition 2.7 in Ch. 1 of \cite{CF-AV} we see that $\l_{\tilde U} = \varphi_{\L}$.

To show that $\l_{\tilde U}$ is a polarization we apply Corollary VIII 7 of \cite{R-AShGrSch} as in the proof of the preceding lemma.
\end{proof}
\section{Kuga-Satake morphisms over fields of characteristic zero}
In the following, sections we will recall the construction of Kuga-Satake abelian varieties associated to polarized K3 surfaces. In our exposition we will follow \cite{D-K3} and \cite{A-HV}. In Section \ref{KSMorphismC} we will use these ideas and the results of \cite{JR-CMK3} to define Kuga-Satake morphisms over number fields.

\subsection{Clifford groups}
Clifford groups will play an essential role in the construction of the Kuga-Satake morphisms and in this section we give a short review of the results we use later in this note. For details we refer to Ch. V in~\cite{Lam-QF} and Ch. 9 in~\cite{WS-QHF}.

Let $d$ be a natural number. For simplification we change the notations of Section 2.1 in \cite{Riz-MK3} by setting $(L,\psi)$ to be the lattice $(L_{2d},\psi_{2d})$ and $(V,\psi)$ to be the quadratic space $(L_{2d},\psi_{2d})\otimes \Q$.

Denote by $G$ the algebraic group $\SO(V,\psi) \cong \SO(2,19)$ over $\Q$ (see Section 2.1 in \cite{JR-CMK3}). Further, following the notations of Example 5.1.4 in \cite{Riz-MK3} we consider the \emph{even Clifford algebra} $C^+(V)$ and the \emph{Clifford group} $G_1 := \CSpin(V)$ of $(V, \psi)$. Recall that one has a homomorphism of algebraic groups
\begin{equation}\label{AdjRep}
 \a \colon \CSpin(V) \lr \SO(V,\psi)
\end{equation}
defined by
\begin{displaymath}
 \a(g) = (v \mapsto gvg^{-1})
\end{displaymath}
which is called the \emph{adjoint representation} of $\CSpin(V)$ on $V$. The kernel of the adjoint representation of $\CSpin(V)$ is $\Gm_m$ and one has a short exact sequence
\begin{equation}\label{CSpinSO}
 1 \lr \Gm_m \lr G_1 \lr G \lr 1.
\end{equation}
Then $G = G_1^{\rm ad}$ and we have that the center $Z(G_1)$ of $G_1$ is $\Gm_m$. 

There is a canonical involution 
\begin{displaymath}
 \iota \colon C^+(L) \lr C^+(L)
\end{displaymath}
which acts trivially on the constants $\Gm_m \hookrightarrow G_1$. We define the \emph{spinorial norm}
\begin{equation}\label{SpNomrm}
\text{N} \colon G_1 \lr \Gm_m
\end{equation}
by setting
\begin{displaymath}
 \text{N}(g) = \iota(g)g.
\end{displaymath} 
It is a surjective homomorphism and we denote its kernel by $\text{Spin}(V)$. The spinorial norm gives rise to a short exact sequence
\begin{displaymath}
1 \lr \text{Spin}(V) \lr G_1 \lr \Gm_m \lr 1.
\end{displaymath} 
One has that $G_1^{\rm der} = \text{Spin}(V)$ is the derived group of $G$ and $\text{N}\colon G_1 \lr \Gm_m = G^{\rm ab}_1$ is the maximal abelian quotient of $G_1$. The group $\text{Spin}(V)$ is simply connected.

Further, we dispose of an embedding $\CSpin(V) \hookrightarrow C^+(V)^*$ and left multiplication by elements of $\CSpin(V)$ on $C^+(V)$ gives an inclusion of algebraic groups 
\begin{equation}\label{SpinRep}
\beta \colon \CSpin(V) \hookrightarrow \GL(C^+(V)).
\end{equation}
See \S 3.2 in \cite{D-K3}. It is called the \emph{spin representation} of $\CSpin(V)$ on $C^+(V)$.
\subsection{Kuga-Satake abelian varieties associated to polarized K3 surfaces}\label{KSPpsect}
In this section we recall the construction of Kuga-Satake abelian varieties. We will follow closely \cite{D-K3} and \cite{A-HV}.

Let $(X,\L)$ be a primitively polarized complex K3 surface of degree $2d$. We use the notations established in Section 2.1 in \cite{Riz-MK3}. As pointed out in Remark 2.3.2 in {\it loc. cite.}, we can find a marking $m \colon H^2_B(X,\Z(1)) \lr L_0$ such that $m(c_1(\L)) = e_1 - df_1$. Then, we obtain an isometry $m \colon P^2_B(X,\Z(1)) \lr L$ and hence the homomorphism $h_X \colon \s \lr \SO\bigl(P^2_B(X,\Z(1))\bigr)$ defines an element
\begin{displaymath}
 h_m : = m \circ h_X \circ m^{-1} \colon \s \lr \SO(V_\R)
\end{displaymath}
of $\Omega^\pm$. There is a unique homomorphism
\begin{displaymath}
 \tilde h_m \colon \s \lr G_{1,\R}
\end{displaymath}
such that $h_m = \a \circ \tilde h_m$, where $\a \colon G_1\lr G$ is the adjoint representation homomorphism (see \S 4.2 in \cite{D-K3}). Let $W$ denote the $\Z$-module $C^+(L)$. The composition of the homomorphism $\tilde h_m$ with the spin representation $\beta \colon G_1 \hookrightarrow \GL(W_\R)$ 
$$
 \beta \circ \tilde h_m \colon \s \lr \GL(W_\R)
$$ 
gives rise to a polarizable $\Z$-HS of type $\{(1,0),(0,1)\}$ on $W$. We refer to Proposition 4.5 in \cite{D-K3} for a proof. Hence $\beta \circ \tilde h_m$ defines a complex abelian variety $A = A(L, h)$, given by the condition that $H^1_B(A,\Z) = W$ as $\Z$-HS. Its dimension is $g = 2^{19}$. 

If we take a different marking $m' \colon H^2_B(X,\Z(1)) \lr L_0$ with $m'(c_1(\l)) = e_1-df_1$, then we have that $m'\circ m^{-1} \circ h_m(z) = h_{m'}(z) \circ m'\circ m^{-1}$ for all $z \in \s$. Therefore $C^+(m'\circ m) \colon W \lr W$ defines an isomorphism between the $\Z$-HS on $W$ induced by $\beta \circ \tilde h_m$ and $\beta \circ \tilde h_{m'}$. Hence we obtain an isomorphism between the abelian varieties associated to $(W,\beta \circ \tilde h_m)$ and $(W,\beta \circ \tilde h_{m'})$. Thus we see that the construction described above associates to a polarized K3 surface $(X,\L)$ an abelian variety $A$, which does not depend on the choice of a marking $m$.
\begin{dfn}
The abelian variety $A$ is called the \emph{Kuga-Satake abelian variety} associated to $(X, \L)$.
\end{dfn}
We will see in Section \ref{PolKSAss} how to give explicitly polarizations of $A$.
\begin{exa}
We shall describe explicitly how to obtain the Hodge structure on $C^+(V)$ in terms of the one on $V$. Choose an orthonormal basis $(e_1,e_2)$ of $V_+ = V_\R \cap (V^{-1,1}\oplus V^{1,-1})$ and let $e_+ = e_1e_2$. Choose an orientation of $(e_1,e_2)$ such that $e_1-ie_2$ spans $V^{1,-1}$. Then multiplication by $e_+$
$$
 x \mapsto e_+x
$$
defines a complex structure on $C^+(V_\R)$ which corresponds to the morphism $\tilde h \colon \s \lr \GL(C^+(V_\R))$ defined above. The Kuga-Satake abelian variety $A$ associated to $(X,\L)$ is exactly the complex torus $C^+(V_\R)/C^+(L)$ where $C^+(V_\R)$ is considered as a complex vector space via the complex structure given by multiplication by $e_+$. For further details we refer to the articles of Satake \cite{S-KS}, Kuga and Satake \cite{KS-KS} and van Geemen \cite{BvG}.
\end{exa}
\begin{exa}
As before $X$ will be a complex K3 surface. Instead of taking the orthogonal complement of an ample line bundle one can consider a subgroup $N \subset c_1(\Pic(X)) \subset H^2_B(X,\Z(1))$ and its complement 
$$
 L_N = N^\perp \subset H^2_B(X,\Z(1))
$$ 
with respect to the bilinear form $\psi$. Then $L_N$ carries a natural polarized $\Z$-HS of type $\{(1,-1),(0,0),(-1,1)\}$ and one can consider again $C^+(L_N)$ and give it a polarized $\Z$-HS of type $\{(1,0),(0,1)\}$ as above. It gives rise to a complex abelian variety $A_N$ associated to the pair $(X,N)$. We refer to \S 4.1 in \cite{M-KSKum} for further comments.

For two subgroups $N \subset N' \subset \NS(X)$ with $d = \dim_\Q({N'}_\Q/N_\Q)$ one has that $A_N$ is isogenous to a product of $2^d$ copies of $A_{N'}$. For a proof see \S 4.4 in \cite{M-KSKum}.
\end{exa}
\begin{exa}
Let $X$ be an exceptional K3 surface. Then the transcendental space $T_{X,\Q} = c_1(\Pic(X)_\Q)^\perp$ is of dimension 2 over $\Q$. By the preceding remarks we conclude that $A$ is isogenous to a product of $2^{19}$ copies of an elliptic curve $E$ which has complex multiplication. See also pp. 241-242 in \cite{KS-KS}. 
\end{exa}
\begin{rem}
Note that from the very construction of $A$ we have that the Mumford-Tate group $\MT(A)$ is contained in $G_1$ viewed as a subgroup of $\GL(C^+(V))$ via the spin representation \eqref{SpinRep}. Moreover, from the short exact sequence (\ref{CSpinSO}) we see that $\Gm_m = \ker(\a)$ is contained in $\MT(A)$ for weight reasons and hence we have an exact sequence
$$
 1 \lr \Gm_m \lr \MT(A) \lr \MT(X) \lr 1.
$$
We also conclude that $\Hg(A)$ is an extension of $\Hg(X)$ by $\Z/2\Z$. 
\end{rem}
\subsection{Endomorphisms}
Denote by $C^+$ the opposite ring $C^+(L)^{\rm op}$. It is non-canonically isomorphic to $C^+(L)$. Let $(X,\L)$ be a primitively polarized K3 surface. Fix a marking $m \colon P^2_B(X,\Z(1)) \lr L$ as in Section \ref{KSPpsect} and let $h_m := m \circ h_X \circ m^{-1} \colon \s \lr G_\R$ be the homomorphism defining the $\Z$-HS on $L$. The right action of $C^+$ on $W:= C^+(L)$ commutes with the morphism $\beta \circ \tilde h_m$ so the Kuga-Satake abelian variety $A$ has complex multiplication by $C^+$ (cf. \S 4.2 in \cite{A-HV} and \S 3.3 in \cite{D-K3}). In other words there is an injection
\begin{equation}\label{EndKSAV}
 \gamma \colon C^+ \hookrightarrow \End(A).
\end{equation}
In fact one can see that there is an isomorphism of $\Z$-HS of type $\{(-1,1),(0,0),(1,-1)\}$
$$
 \phi_\Z \colon C^+(L)_{\rm ad} \lr \End_{C^+}(W)
$$
where $C^+(L)_{\rm ad}$ is the $\Z$-HS obtained from $(L,h)$ using the tensor construction $C^+(\ )$.
\subsection{Polarizations}\label{PolKSAss}
Let $(X,\L)$ be a complex K3 surface with a primitive polarization $\l$ of degree $2d$ and let $m\colon H^2_B(X,\Z(1)) \lr L_0$ be a marking such that $m(c_1(\L)) = e_1-df_1$ (see Section 2.1 and Remark 2.3.2 in \cite{Riz-MK3}). Let $h_m \colon \s \lr \SO(V_\R)$ be the $\Z$-HS induced on $L$ by $h_X$ and let let $A$ be its associated Kuga-Satake abelian variety. We will show how to give explicitly a polarization of the $\Z$-HS $W$ ($=C^+(L) = H^1_B(A,\Z)$). 

Let $\iota \colon C^+(L) \lr C^+(L)$ be the canonical involution of the even Clifford algebra. Fix a non-zero element $a \in C^+$ such that $\iota(a) = -a$. Then the skew-symmetric form 
\begin{equation}\label{PolFormKS}
 \varphi_a \colon W \otimes W \lr \Z(-1)
\end{equation} 
given by 
\begin{displaymath}
 \varphi_a(x,y) = \text{tr}(\iota(x)ya)
\end{displaymath}
defines a polarization for the $\Z$-HS on $W$ if and only if the symmetric bilinear form $i\varphi_a(x,\tilde h_m(i)y)$ is positive definite (here $i = \sqrt {-1}$). Lemma 4.3 in \cite{D-K3} (see also Example \ref{PolExample}) guarantees the existence of an element $a \in C^+$ for which $\pm \varphi_a$ is a polarization.
\begin{exa}\label{PolExample}
Let $e_1, \dots, e_{21}$ be an orthogonal basis of $(V, \varphi)$ such that $\psi(e_i, e_i) < 0$ for $i = 1,2$. Let $m \ne 0$ be an integer such that $me_1e_2 \in C^+(L)$. One has that $\iota(me_1e_2) = -me_1e_2$ and if $h \in \Omega^{\pm}$, then either $\varphi_{me_1e_2}$ or $-\varphi_{me_1e_2}$ is a polarization for $\tilde h$. For a proof we refer to \cite[Prop. 5.9]{BvG}.
\end{exa}
\begin{rem}\label{PolRem}
Note that the degree of the polarization $\varphi_a$ depends only on $a$ and $d$ and can be computed explicitly.
\end{rem}
\begin{rem}
Let $a \in C^+$ be an element such that $\iota(a) = -a$ and, say $\varphi_a$ is a polarization the $\Z$-HS on $W$ induced by $\tilde h$. Then $\varphi_a$ defines a polarization $\mu \colon A \lr A^t$ which gives rise to a Rosati involution $\dagger$ on $\End^0(A) = \End(A) \otimes \Q$. One can see that the restriction of the Rosati involution to $C^+\otimes \Q \hookrightarrow \End^0(A)$ (cf. \eqref{EndKSAV}) is given by
\begin{displaymath}
f^\dagger = a^{-1} \iota(f) a
\end{displaymath} 
for all $f \in C^+(V)$. Hence $C^+(V)$ is stable under $\dagger$.
\end{rem}
\begin{rem}\label{MarkingDepPol}
Note that we make some non-canonical choices to define a polarization on $A$. For instance, it is not clear if two different markings $m_i \colon H^2_B(X,\Z(1)) \lr L_0$ for which $m_i(c_1(\l)) = e_1 - df_1$ give rise to two isomorphic \emph{polarized} abelian varieties $(A,\mu_1)$ and $(A,\mu_2)$. 
\end{rem}
\subsection{Kuga-Satake morphisms over fields of characteristic zero}\label{KSMorphismC}
Recall that we associated to every polarized complex K3 surface a complex abelian variety. We will explain here how to do this in families. Following the line of thoughts in \cite{D-K3} and \cite{A-HV} we define Kuga-Satake morphisms from the moduli spaces of polarized K3 surfaces with certain level structures to moduli stacks of polarized abelian varieties. We shall keep the notations from the previous sections.

Consider the Shimura datum $(G,\Omega^\pm)$ (cf. Section 2.2 in \cite{JR-CMK3}) and let $h_0 \colon \s \lr G_R$ be an element of $\Omega^\pm$. Let $\tilde h_0 \colon \s \lr G_{1,\R}$ be the unique homomorphism such that $h_0 = \a \circ \tilde h_0$ (cf. Section \ref{KSPpsect}). Define $\Omega^\pm_1$ to be the $G_1(\R)$-conjugacy class of $\tilde h_0$. The pair $(G_1,\Omega^\pm_1)$ defines a Shimura datum with reflex field $\Q$. We refer to Appendix 1 in \cite{A-HV} for a proof. 

The adjoint representation \eqref{AdjRep} defines a morphism of Shimura data
$$
 \a \colon (G_1,\Omega_1^\pm) \lr (G,\Omega^\pm)
$$
in the following way: $\a_{gr} \colon G_1 \lr G$ is the adjoint representation homomorphism and $\a_{HS} \colon \Omega^\pm_1 \lr \Omega^\pm$ the the  morphism sending $\tilde h$ to $\a \circ \tilde h$. The morphism $\a_{HS}$ is well-defined as $h = g\circ \tilde h_0 \circ g^{-1}$ for some $g \in G_1(\R)$ and hence $\a \circ \tilde h = \a(g) \circ h_0 \circ \a(g)^{-1} \in \Omega^\pm$. Moreover $\a_{HS} \colon \Omega^\pm_1 \lr \Omega^\pm$ is an analytic isomorphism (\S 4.2 in \cite{D-K3} or Lemma 4.11 in \cite{Mil-ShVGRed}).

We use the notations established in Section 5 of \cite{Riz-MK3}. Fix a natural number $n \geq 3$. Let $\Kg^{\rm sp} \subset G_1(\A_f)$ be a subgroup of finite index in $\Kg^{\rm sp}_n$ and denote by $\Kg^{\rm a}$ the image $\a(\Kg^{\rm sp}) \subset G(\A_f)$ which is a subgroup of finite index in $\Kg^{\rm a}_n$ (cf. Example 5.1.4 in \cite{Riz-MK3}). Then one has a morphism of quasi-projective $\Q$-schemes
\begin{equation}\label{SpOrthMap}
 \a_{(\Kg^{\rm sp}, \Kg^{\rm a})} \colon Sh_{\Kg^{\rm sp}}(G_1, \Omega^{\pm}_1)  \lr Sh_{\Kg^{\rm a}}(G,\Omega^\pm).
\end{equation}
Consider the group $C = \Gm_m(\Q)\backslash \a^{-1}(\Kg^{\rm a})/\Kg^{\rm sp}$. We have that 
\begin{equation}\label{ShGalGroiup}
\begin{split}
C &= \Gm_m(\Q)\backslash \a^{-1}(\Kg^{\rm a})/\Kg^{\rm sp} = \Gm_m(\Q)\backslash \Gm_m(\A_f)\Kg^{\rm sp}/\Kg^{\rm sp} \\
  &= \Gm_m(\Q)\backslash \Gm_m(\Q)\Gm_m(\hat \Z)\Kg^{\rm sp}/\Kg^{\rm sp} \cong \Gm_m(\hat \Z)/(\Gm_m(\hat \Z) \cap \Kg^{\rm sp}).
\end{split}
\end{equation}
The group $C$ acts on $Sh_{\Kg^{\rm sp}}(G_1,\Omega^\pm_1)_\C$ via right multiplication. We have that $Z(G_1) = \Gm_m$ and $G = G_1/Z(G_1)$. Further, by Hilbert's Theorem 90, $H^1(k,\Gm_m) = 0$ for all fields of characteristic zero, hence we can apply Lemma 4.13 in \cite{Mil-ShVGRed} and conclude that the morphism $\a_{(\Kg^{\rm sp}, \Kg^{\rm a})} \otimes \C$ is a Galois cover with a Galois group $C$. As $C$ acts on $Sh_{\Kg^{\rm sp}}(G_1,\Omega^\pm_1)_\C$ via Hecke correspondences we see that these automorphisms are defined over $\Q$. Therefore the morphism \eqref{SpOrthMap} is a Galois cover with a Galois group $C$.

We will describe more explicitly the relation between these two Shimura varieties over $\C$. 
Consider the finite sets $\mathcal C_{G_1} := G_1(\Q)_+\backslash G_1(\A_f)/\Kg^{\rm sp}$ and $\mathcal C_G := G(\Q)_+\backslash G(\A_f)/\Kg^{\rm a}$. The homomorphism $\a$ defines a surjective map of sets $\a \colon \mathcal C_{G_1} \lr \mathcal C_G$ (cf. \eqref{CSpinSO}). Note that $C$ naturally acts on $\mathcal C_{G_1}$ from the right. With this action the map $\a$ makes $\mathcal C_{G_1}$ into a $C$-torsor over $\mathcal C_{G}$ (in the sense of sets); in other words, if $[g] \in \mathcal C_G$ and $g_1 \in G_1(\A_f)$ is an element with $\a([g_1]) = [g]$, then the map $C \lr \a^{-1}([g])$ given by $u \mapsto [g_1u]$ is a bijection.

One has that the decomposition of $Sh_{\Kg^{\rm sp}}(G_1,\Omega_1^\pm)_\C$ into connected components is
\begin{displaymath}
Sh_{\Kg^{\rm sp}}(G_1,\Omega^\pm_1)_\C = \coprod_{[g] \in \mathcal C_{G_1}} {\Gamma'}_{[g]} \backslash \Omega_1^+
\end{displaymath}
where ${\Gamma'}_{[g]}= G_1(\Q)_+ \cap g\Kg^{\rm sp}g^{-1}$, for some representative $g$ of the class $[g]$ (see \cite[\S 5, Lemma 5.13]{Mil-IShV}. Similarly, we have that
\begin{displaymath}
Sh_{\Kg^{\rm a}}(G,\Omega^\pm)_\C = \coprod_{[g] \in \mathcal C_G} \Gamma_{[g]} \backslash \Omega^+
\end{displaymath}
where $\Gamma_{[g]}= G(\Q)_+ \cap g\Kg^{\rm a}g^{-1}$ for some representative $g$ of $[g]$. 

The morphism $\a_{(\Kg^{\rm sp},\Kg^{\rm a})}$ maps the connected component ${\Gamma'}_{[g]} \backslash \Omega_1^+$ to ${\Gamma}_{[\a(g)]} \backslash \Omega^+$ sending the class $[\tilde h]$ to the class $[h]$ (cf. \S \ref{KSPpsect}). The restriction
\begin{equation}\label{RestSpOrthMap}
\a_{(\Kg^{\rm sp},\Kg^{\rm a})} \colon  {\Gamma'}_{[g]} \backslash \Omega_1^+ \lr {\Gamma}_{[\a(g)]} \backslash \Omega^+
\end{equation}
is an isomorphism of complex quasi-projective varieties. Indeed, $\a$ maps ${\Gamma'}_{[g]}$ surjectively onto ${\Gamma}_{[\a(g)]}$ and as $-1 \not \in {\Gamma'}_{[g]}$ (because $-1 \not \in g\Kg^{\rm sp}_ng^{-1} \supset g\Kg^{\rm sp}g^{-1}$) one concludes, from the short exact sequence (8) in Example 5.1.4 in \cite{Riz-MK3}, that ${\Gamma'}_{[g]}$ is mapped isomorphically onto $\Gamma_{[\a(g)]}$. The morphism $\a_{HS} \colon \Omega^+_1 \lr \Omega^+$ is an isomorphism so we see that \eqref{RestSpOrthMap} is an isomorphism as well. Further, we have that
\begin{equation}\label{ConnCptsDecomp}
 \a_{(\Kg^{\rm sp},\Kg^{\rm a})}^{-1}(\Gamma_{[g]} \backslash \Omega^+) = \coprod_{u \in C} {\Gamma'}_{[g_1u]} \backslash \Omega_1^+
\end{equation}
where $g_1 \in G(\A_f)$ with $\a(g_1) = g$.

Denote by $W$ the $\Z$-module $C^+(L)$ and choose an element $a \in C^+$ such that $\iota(a) = -a$. Recall that for such an element we have defined a bilinear form $\varphi_a \colon W \otimes W \lr \Z(-1)$ (see \eqref{PolFormKS}). The image of $G_1$ under the spin representation $\beta \colon G_1 \hookrightarrow \GL(W_\Q)$ is actually contained in $\CSp(W_\Q,\varphi_a)$. Indeed, for any element $\gamma \in G_1$ we have that
\begin{displaymath}
\begin{split}
\varphi_a(\gamma x,\gamma y) &= \text{tr}(\iota(\gamma x)\gamma y a) = \text{tr}(\iota(x)\iota(\gamma) \gamma y a) \\
&= \text{tr}\bigl(\iota(x){\rm N}(\gamma)y a\bigr) = {\rm N}(\gamma)\text{tr}(\iota(x)ya) \\
&= {\rm N}(\gamma)\varphi_a(x,y)\\
\end{split}
\end{displaymath}
hence $\beta(\gamma) \in \CSp(W_\Q,\varphi_a)$. Further, if the bilinear form $
\varphi_a$ defines a polarization for a Hodge structure $\beta \circ \tilde h_1$ on $W$, then it defines a polarization for all Hodge structures $\beta \circ \tilde h$ on $W$, for which $\tilde h$ belongs to the connected component of $\Omega^\pm_1$ of $\tilde h$. If $\varphi_a$ is a polarization for those $\tilde h$ coming from the elements in $\Omega^+$, then $-\varphi_a$ is a polarization for the $\tilde h$ coming from the elements in $\Omega^-$.
\begin{assmt}\label{PolAssumpt}
We assume that $a \in C^+$ is such that $\iota(a) = -a$ and that $\varphi_a$ or $-\varphi_a$ defines a polarization for the $\Z$-Hodge structures induced on $W$ by $\beta \circ \tilde h$ for any $\tilde h \in \Omega^\pm_1$. 
\end{assmt}
Under the above assumption define the inclusion of Shimura data
$$
 \beta \colon (G_1,\Omega^\pm) \hookrightarrow (\CSp(W_\Q,\varphi_a), \mathfrak{H}^{\pm}).
$$
as $\beta_{gr} \colon G_1 \hookrightarrow \CSp(W_\Q,\varphi_a)$ being the spin representation \eqref{SpinRep} and $\beta_{HS} \colon \Omega^\pm_1 \hookrightarrow \mathfrak{H}^{\pm}$ mapping $\tilde h$ to $\beta \circ \tilde h$.

Let $\Lambda_n$ be the congruence level $n$-subgroup of $\CSp(W_\Q,\varphi_a)(\A_f)$ corresponding to the lattice $W$ of $W_\Q$. In other words we take
\begin{displaymath}
 \Lambda_n = \{ g\in \CSp(W_\Q,\varphi_a)(\A_f)\ |\ gW_{\hat \Z}=W_{\hat \Z}\ \text{and}\ g \equiv 1 \pmod n\}.
\end{displaymath} 
It is clear from the definitions that $\beta(\Kg^{\rm sp}) \subset  \beta(\Kg^{\rm sp}_n) \subset \Lambda_n$ hence we obtain a morphism of quasi-projective $\Q$-schemes  
\begin{displaymath}
\beta_{(\Kg^{\rm sp}, \Lambda_n)} \colon Sh_{\Kg^{\rm sp}}(G_1, \Omega^{\pm}) \lr Sh_{\Kg^{\rm sp}_n}(G_1, \Omega^{\pm}) \lr Sh_{\Lambda_n}(\CSp(W_\Q,\varphi_a), \mathfrak{H}^{\pm}).
\end{displaymath}
Note that fixing the lattice $W$, respectively the arithmetic group $\Lambda_n$, one has an immersion $Sh_{\Lambda_n}(\CSp(W_\Q,\varphi_a), \mathfrak{H}^\pm) \hookrightarrow \Av_{g,d',n,\Q}$ where $d'$ is explicitly computed in terms of $d$ and $a$ (cf. Remark \ref{PolRem}). It is given by the identification of $Sh_{\Lambda_n}(\CSp(W_\Q,\varphi_a), \mathfrak{H}^\pm)$ with a component $\Av_{g,\delta,n,\Q}$ of $\Av_{g,d',n,\Q}$ corresponding to an elementary divisor sequence $\delta = (d_1,\dots,d_r)$, uniquely determined by $\varphi_a$, with $d_1\cdots d_r = d'$ (cf. Definition 1.3 in \cite{deJ-AV}). We can put all morphisms considered so far in the following diagram
\begin{equation}\label{KSDefDiagram}
\xymatrix{ 
 Sh_{\Kg^{\rm sp}}(G_1,\Omega_1^\pm) \ar[d]^{\a_{(\Kg^{\rm sp},\Kg^{\rm a})}} \ar[rrr]^{\beta_{(\Kg^{\rm sp},\Lambda_n)}} & & & \Av_{g,d',n,\Q} \ar[d]^{pr_n} \\
 Sh_{\Kg^{\rm a}}(G,\Omega^\pm) & & & \Av_{g,d',\Q}.
}
\end{equation}
\newline
\newline
{\bf First construction of Kuga-Satake morphisms.} Both morphisms $\a_{(\Kg^{\rm sp}, \Kg^{\rm a})}$ and $pr_n$ are quotient morphisms as $\Av_{g,d',\Z[1/n]}$ is the quotient stack $[\Av_{g,d',n}/\GL(W_\Q)(\Z/n\Z)]$ (cf. \S 4.3.4 in Ch VII of \cite{MB-AbV}). Moreover, as $C$ acts freely on $Sh_{\Kg^{\rm sp}}(G_1,\Omega^\pm)$ we have that the stack $[Sh_{\Kg^{\rm sp}}(G_1,\Omega_1^\pm)/C]$ is represented by the quotient scheme $Sh_{\Kg^{\rm a}}(G_1,\Omega^\pm) \cong Sh_{\Kg^{\rm sp}}(G_1,\Omega_1^\pm)/C$. The spin representation defines a homomorphism (see \eqref{ShGalGroiup})
\begin{equation}\label{GalGrHom}
\beta \colon C = \Gm_m(\hat \Z)\Kg^{\rm sp}/\Kg^{\rm sp} \lr \GL(W_\Q)(\Z/n\Z).
\end{equation} 
We will show that $\beta_{(\Kg^{\rm sp}, \Lambda_n)}$ descends to a morphism $\beta_{\Kg^{\rm a}} \colon Sh_{\Kg^{\rm a}}(G,\Omega^\pm) \lr \Av_{g,d',\Q}$. To do this we have to check that $\beta_{(\Kg^{\rm sp}, \Lambda_n)}$ is equivariant with respect to the homomorphism \eqref{GalGrHom}. Both $Sh_{\Kg^{\rm a}}(G,\Omega^\pm)$ and $\Av_{g,d',n,\Q}$ are reduced schemes over $\Q$ so we can check the statement on $\C$-valued points. In other words we have to show that
\begin{displaymath}
 \beta_{(\Kg^{\rm sp},\Lambda_n)}\bigl(g \cdot [\tilde h, r]_{\Kg^{\rm sp}}\bigr) = \beta(g)\cdot \beta_{(\Kg^{\rm sp},\Lambda_n)}\bigl([\tilde h, r]_{\Kg^{\rm sp}}\bigr)
\end{displaymath}
for any $g \in C, \tilde h \in \Omega^\pm$ and $r \in G_1(\A_f)$. But this is tautology as from the definitions we see that
\begin{displaymath}
\begin{split}
\beta_{(\Kg^{\rm sp},\Lambda_n)}\bigl(g \cdot [\tilde h, r]_{\Kg^{\rm sp}}\bigr) &=  \beta_{(\Kg^{\rm sp},\Lambda_n)}\bigl([\tilde h, rg]_{\Kg^{\rm sp}}\bigr) = [\beta \circ \tilde h, \beta(rg)]_{\Lambda_n} \\
 &= \beta(g)\cdot [\beta \circ \tilde h, \beta(r)]_{\Lambda_n}= \beta(g) \cdot \beta_{(\Kg^{\rm sp},\Lambda_n)}\bigl([\tilde h, r]_{\Kg^{\rm sp}}\bigr).\\
\end{split}
\end{displaymath}
Hence $\beta_{(\Kg^{\rm sp}, \Lambda_n)}$ descends to a morphism of algebraic stacks
\begin{displaymath}
 \beta_{\Kg^{\rm a}} \colon Sh_{\Kg^{\rm a}}(G,\Omega^\pm) \lr \Av_{g,d',\Q}.
\end{displaymath}

Recall that in Section 2.4 of \cite{JR-CMK3} we defined a period morphism $j_{d,\Kg^{\rm a}} \colon \Fk_{2d,\Kg^{\rm a},\Q} \lr Sh_{\Kg^{\rm a}}(G,\Omega^\pm)$ which sends any complex polarized K3 surface with a level $\Kg^{\rm a}$-structure to its period point.
\begin{dfn}
Define the \emph{Kuga-Satake morphism} associated to $d,a$ and $\Kg^{\rm a}$
\begin{displaymath}
 f^{ks}_{d,a,\Kg^{\rm a},\Q} \colon \Fk_{2d,\Kg^{\rm a},\Q} \lr \Av_{g,d',\Q}
\end{displaymath}
to be the composite $f^{ks}_{d,a,\Kg^{\rm a},\Q} = \beta_{\Kg^{\rm a}} \circ j_{d,\Kg^{\rm a}}$.
\end{dfn}
Thus we have proved the following statement. 
\begin{prp}\label{KSCmap}
Let $d, n\in \N$ with $n \geq 3$ and let $\Kg^{{\rm sp}} \subset \Kg^{\rm sp}_n$ be a subgroup of finite index. Fix a non-zero element $a \in C^+$ which satisfies Assumption \ref{PolAssumpt}. Then one has a Kuga-Satake morphism
$$
 f_{d,a,\Kg^{\rm a}, \Q}^{ks} \colon \Fk_{2d, \Kg^{\rm a}, \Q} \lr \Av_{g,d',\Q}
$$
where $g = 2^{19}$ and $d'$ depends explicitly on $a$ and $d$. It maps every primitively polarized complex K3 surface $(X,\l,\nu)$ with a level $\Kg^{\rm a}$-structure $\nu$ to its associated Kuga-Satake abelian variety $A$ with a certain polarization of degree ${d'}^2$. 
\end{prp}
\begin{rem}
Note that if $\Kg_1$ is a subgroup of $G_1(\hat \Z)$ of finite index contained in $\Kg^{\rm sp}_n$ and such that $\a(\Kg_1) = \Kg^{\rm a}$, then the morphism $\beta_{(\Kg_1,\Lambda_n)} \colon Sh_{\Kg_1}(G_1,\Omega_1^\pm) \lr \Av_{g,d',n,\Q}$ also descends to the morphism $\beta_{\Kg^{\rm a}} \colon Sh_{\Kg^{\rm a}}(G,\Omega^\pm) \lr \Av_{g,d',\Q}$. 
\end{rem}
\begin{exa}
Take $\Kg^{\rm sp}$ to be the group $\Kg^{\rm sp}_n$. Then $\Kg^{\rm a} = \Kg^{\rm a}_n$ and we obtain a Kuga-Satake morphism
\begin{displaymath}
f^{ks}_{d,a,n,\Q} \colon \Fk_{2d,n^{\rm sp},\Q} \lr \Av_{g,d',\Q}.
\end{displaymath}
\end{exa}
\begin{rem}
If $\Kg^{\rm a}$ is an admissible subgroup of $\SO(V)(\hat \Z)$ (see Definition 5.3.1 in \cite{Riz-MK3}), then we have an open immersion $j_{d,\Kg^{\rm a}} \colon \Fk^{\rm full}_{2d,\Kg^{\rm a},\Q} \hookrightarrow Sh_{\Kg^{\rm a}}(G,\Omega^\pm)$ and therefore we obtain a Kuga-Satake morphism
$$
 f^{ks}_{d,a,\Kg^{\rm a},\Q} \colon \Fk^{\rm full}_{2d,\Kg^{\rm a},\Q} \lr \Av_{g,d',\Q}
$$
defined by $f^{ks}_{d,a,\Kg^{\rm a},\Q} = \beta_{\Kg^{\rm a}} \circ j_{d,\Kg^{\rm a}}$.
\end{rem}
\begin{rem}
One might want to descend the Kuga-Satake morphism defined in Proposition \ref{KSCmap} to a morphism $\Fk_{2d,\Q} \lr \Av_{g,d',\Q}$. The essence of the problem is that the Kuga-Satake construction described above requires a non-canonical choice of an element $a \in C^+$ to define a polarization. One can show that the obstruction for descending the Kuga-Satake morphism to a map $\Fk_{2d,\Q} \lr \Av_{g,d',\Q}$ is equivalent to the problem posed in Remark \ref{MarkingDepPol}.
\end{rem}
Our main goal is to define Kuga-Satake morphisms in mixed characteristic. As we will see later (Remark \ref{ProblemExtKSQ}) there are problems extending the morphism $f^{ks}_{d,a,n} \colon \Fk_{2d,n^{\rm sp},\Q} \lr \Av_{g,d',\Q}$ due to the fact that $\Av_{g,d'}$ is an algebraic stack. We will give a second construction of Kuga-Satake morphisms below to which we can apply the extension result of Section \ref{AbstractExtension}.
\newline
\newline
{\bf Second construction of Kuga-Satake morphisms.} We will construct a morphism $f^{ks}_{d,a,\gamma,n,E} \colon \Fk_{2d,n^{\rm sp},E} \lr \Av_{g,d',n,E}$ for a number field $E$ which can be determined via class field theory from the data $d,a,\gamma,n$ (see below). To do that we will first determine the fields of definition of the geometric connected components of $Sh_{\Kg^{\rm sp}_n}(G_1,\Omega^\pm_1)$ and  $Sh_{\Kg^{\rm a}_n}(G,\Omega^\pm)$.

We have that
\begin{displaymath}
 \pi_0\bigl(Sh_{\Kg^{\rm sp}_n}(G_1,\Omega^\pm_1)_\C\bigr) \cong G_1(\Q)_+\backslash G_1(\A_f)/\Kg^{\rm sp}_n \cong \Gm_m(\A)/\bigl(\Q^\times \R_{>0} \text{N}(\Kg^{\rm sp}_n)\bigr)
\end{displaymath}
where $\text{N} \colon G_1 \lr G_1^{\rm ab} = \Gm_m$ is the spinorial norm homomorphism (see \eqref{SpNomrm}). Denote by $E_n$ the subfield of $\Q^{\rm ab}$ corresponding to the group $\Q^\times \R_{>0}\text{N}(\Kg^{\rm sp}_n)$ via class field theory (see Section 3.1 in \cite{JR-CMK3}). Then, we have an isomorphism
\begin{displaymath}
\text{art}_{E_n/\Q} \colon \Gm_m(\A)/\bigl(\Q^\times \R_{>0}\text{N}(\Kg^{\rm sp}_n)\bigr) \lr \Gal(E_n/\Q).
\end{displaymath}
The Galois action on the geometric connected components of $Sh_{\Kg^{\rm sp}_n}(G_1,\Omega_1^\pm)$ is given as follows: Let $Y$ be the connected component $\Omega_1^+/{\Gamma'}_{[1]}$. It is defined over $E_n$ and if $\sigma \in \Gal(E_n/\Q)$ is an automorphism such that $\text{art}_{E_n/\Q}(\sigma) = \text{N}(g)$ for some $g\in G_1(\A_f)$, then
$Y^\sigma = \Omega_1^+/{\Gamma'}_{[g]}$. We have that
\begin{equation}\label{GalDecCsp}
Sh_{\Kg^{\rm sp}_n}(G_1,\Omega_1^\pm) = \bigcup_{\sigma \in \Gal(E_n/\Q)} Y^\sigma.
\end{equation}
For details see \S 2 in \cite{Kud-ACOShV}. Further, if we denote by $X$ the connected component $\Omega^+/\Gamma_{[1]}$ of $Sh_{\Kg^{\rm a}_n}(G,\Omega^\pm)_\C$, then its field of definition $E_X$ is a subfield of $E_n$ and hence it is an abelian extension of $\Q$. We have that $[E_n:E_X] = \# C = \varphi(n)$ (cf. \eqref{ShGalGroiup}) and
\begin{equation}\label{GalDecSO}
Sh_{\Kg^{\rm a}_n}(G,\Omega^\pm) = \bigcup_{\sigma \in \Gal(E_X/\Q)} X^\sigma.
\end{equation}
In order to define a Kuga-Satake morphism $f^{ks}_{d,a,n,\gamma,E} \colon \Fk_{2d,n^{\rm sp},E} \lr \Av_{g,d',n,E}$ we will give a section of $\a_{(\Kg^{\rm sp}_n,\Kg^{\rm a}_n)}$ and use \eqref{KSDefDiagram}. We see from \eqref{RestSpOrthMap}, \eqref{ConnCptsDecomp}, \eqref{GalDecCsp} and \eqref{GalDecSO}, that giving such a section is equivalent to giving a set-theoretic section of the homomorphism $\Gal(E_n/\Q) \lr \Gal(E_X/\Q)$. For any such (set-theoretic) section $\gamma \colon \Gal(E_X/\Q) \lr \Gal(E_n/\Q)$ one has a morphism
\begin{equation}
 \delta_\gamma \colon Sh_{\Kg^{\rm a}_n}(G,\Omega^\pm) = \bigcup_{\sigma \in \Gal(E_X/\Q)} X^\sigma \cong \bigcup_{\sigma \in \Gal(E_X/\Q)} Y^{\gamma(\sigma)} \subset Sh_{\Kg^{\rm sp}_n}(G_1,\Omega_1^\pm)
\end{equation} 
which is defined over $E_n$.
\begin{dfn}\label{KSLevStrMor}
Define the \emph{Kuga-Satake morphism} associated to $d,a,n$ and $\gamma$
$$
 f^{ks}_{d,a,n,\gamma,E_n} \colon \Fk_{2d,n^{\rm sp},E_n} \lr \Av_{g,d',n,E_n}
$$
to be the composite $f^{ks}_{d,a,n,\gamma,E_n} = \beta_{(\Kg^{\rm sp}_n,\Lambda_n)} \circ \delta_\gamma \circ j_{d,\Kg^{\rm a}_n,E_n}$ defined over $E_n$.
\end{dfn}
Thus, we have proved the following statement. 
\begin{prp}\label{KSCMapN}
Let $d, n\in \N$ with $n \geq 3$. Let $E_n$ and $E_X$ be as above and suppose given a set-theoretic section $\gamma$ of the homomorphism $\Gal(E_n/\Q) \lr \Gal(E_X/\Q)$. Fix a non-zero element $a \in C^+$ which satisfies Assumption \ref{PolAssumpt}. Then one has a Kuga-Satake morphism
$$
 f_{d,a,n,\gamma,E_n}^{ks} \colon \Fk_{2d,n^{\rm sp},E_n} \lr \Av_{g,d',n,E_n}
$$
where $g = 2^{19}$ and $d'$ depends explicitly on $a$ and $d$. It maps every primitively polarized complex K3 surface $(X,\l,\nu)$ with a spin level $n$-structure $\nu$ to its associated Kuga-Satake abelian variety $A$ with a certain polarization of degree ${d'}^2$ and a certain level $n$-structure. Further, by construction, for any choice of a section $\gamma$ we have that $f_{d,a,n}^{ks}\otimes E_n = pr_n \circ f_{d,a,n,\gamma,E_n}^{ks}$.
\end{prp}
As we have seen, there are many possible ways of defining Kuga-Satake morphisms. In general, one has to make some non-canonical choices in order to find a section of $\a_{(\Kg^{\rm sp}_n,\Kg^{\rm a}_n)}$ in \eqref{KSDefDiagram} and define a morphism $\Fk_{2d,n^{sp},\C} \lr \Av_{g,d',n,\C}$. 

Below, we explain the relative Kuga-Satake construction of Deligne (cf. \S 5 in \cite{D-K3}) in our framework. Consider the diagram
\begin{displaymath}
\xymatrix{
\Gamma'_{[1]}\backslash \Omega_1^+ \ar[d]_{\cong}^{\a_{(\Kg^{\rm sp}_n,\Kg^{\rm a}_n)}} \ar[rr]^{\beta_{(\Kg^{\rm sp}_n,\Lambda_n)}} & & \Av_{g,d',n,\C} \\
\Gamma_{[1]} \backslash \Omega^+. &  &
}
\end{displaymath}
Over $\C$ one can define a morphism $\Fk_{2d,n^{sp},\C} \lr \Gamma_{[1]}\backslash \Omega^+$ by mapping all connected components of $\Fk_{2d,n^{sp},\C}$ to $\Gamma_{[1]}\backslash \Omega^+$. See the proof of Proposition 5.7 in \cite{D-K3}. Composing these two maps we obtain a morphism 
$$
 f_n \colon \Fk_{2d,n^{sp},\C} \lr \Gamma_{[1]}\backslash \Omega^+ \cong {\Gamma'}_{[1]}\backslash \Omega^+_1 \lr \Av_{g,d',n,\C}
$$
which is the relative Kuga-Satake construction described in \S 5 in \cite{D-K3} and \S 5 in \cite{A-HV}. One can show that this morphism is defined over a number field. Suppose further, that $n = 3$ or $4$. Here is a possible way to study the field of definition of $f_n$. Combining Proposition 8.3.5 and Theorem 8.4.3 in \cite{A-HV} one can see that $f_n$ is defined over the composite of $E_n$ with any field $K \subset \C$ for which $\Fk_{2d,n^{\rm sp},\Q}$ has a $K$-valued point. Then by Theorem 7 in \cite{Riz-FDef} the morphism $f_n$ is defined over the composite of $E_n$ with the fields of definition of the geometric connected components of $\Fk_{2d,n^{\rm sp},\Q}$. In general, this field can be a non-trivial extension of $E_n$.
\begin{rem}
The construction of Kuga-Satake morphisms described in this section and the one given in \cite{D-K3} and \cite{A-HV} differ in the choice of a period morphism. We use the ``modified'' period map $j_{d,\Kg,\C}$ in order to be able to apply the results of \cite{JR-CMK3}. In this way we can control explicitly the fields of definition of the morphisms involved in the relative Kuga-Satake construction and therefore the field of definition of $ f^{ks}_{d,a,n,\gamma,\C}$.
\end{rem}
We will end this section with a result comparing the \'etale cohomology of a K3 surface and its associated Kuga-Satake abelian variety. Let $U^i$ be a geometric connected component of $\Fk_{2d,n^{sp},E_n}$ which is defined over a field $i \colon K \hookrightarrow \C$. Let $(\pi_{X^i} \colon X^i \lr U^i, \l^i, \nu^i)$ be the pull-back of the universal family to $U^i$. Denote by $(\pi_{A^i} \colon A^i \lr U^i, \mu^i, \e^i)$ the polarized abelian scheme with level $n$-structure $f^{ks}_{d,a,n,\gamma,E_n}((\pi_{X^i} \colon X^i \lr U^i, \l^i, \nu^i))$. 

Taking a base change $i \colon K \lr \C$ we have an abelian scheme $(A^i_\C \lr U^i_\C, \mu^i_\C,\e^i_\C)$ which is exactly $f^{ks}_{d,a,n,\gamma,\C}\bigl((X^i_\C \lr U^i_\C, \l^i_\C, \nu^i_\C)\bigr)$ and which, by construction, has multiplication by $C^+$. Further, we know that $\End_{U^i_\C}(A^i_\C) = C^+$ (see the beginning of \S 8 in \cite{A-HV}) and one has further that $\End_{U^i}(A^i) = C^+$.
\begin{lem}
There is a unique isomorphism of $\Z_l$-sheaves
$$
 C^+(P^2_{\rm et}\pi_{X^i,*}\Z_l(1)) \cong \End_{C^+}(R^1_{\rm et}\pi_{A^i,*}\Z_l).
$$ 
\end{lem}
\begin{proof}
One repeats step by step the proof of Lemma 6.5.13 in \cite{D-K3}.
\end{proof} 
\begin{cor}\label{ChoCompRes}
Let $K$ be a field of characteristic zero and suppose given a $K$-valued point $(X,\l,\nu) \in \Fk_{2d,n^{\rm sp},E_n}(K)$. If $(A,\mu,\e)$ is the corresponding Kuga-Satake abelian variety $f^{ks}_{d,a,n,\gamma,E_n}\bigl((X,\l,\nu)\bigr)$, then one has an isomorphism of $\Gal(\bar K/K)$-modules
$$
 C^+\bigl(P^2_{\rm et}(X_{\bar K}, \Z_l(1))\bigr) \cong \End_{C^+}\bigl(H^1_{\rm et}(A_{\bar K}, \Z_l)\bigr)
$$
for any prime number $l$.
\end{cor}
\begin{proof}
It follows from the preceding lemma.
\end{proof}
\begin{rem}\label{GRedProperty}
Note that if $R$ is a discrete valuation ring with a maximal ideal $\mathfrak p$ and field of fractions $K$ of characteristic zero, containing $E_n$. Suppose given a polarized K3 surface $(X,\l,\a)$ with spin level $n$-structure over $K$ and let $(A,\mu,\b)$ be the corresponding Kuga-Satake abelian variety. Suppose further that $X$ has good reduction modulo $\mathfrak p$. Then the inertia subgroup $I_{\mathfrak p}$ acts trivially on $P^2_{\rm et}(X_{\bar K},\Z_l(1))$ for every $l$ different from the characteristic of $R/\mathfrak p$. As shown in \S 6.6 in \cite{D-K3} and Lemma 9.3.1 in \cite{A-HV} this implies that $I_{\mathfrak p}$ acts via a finite group on $H^1_{\rm et}(A_{\bar K},\Z_l)$ i.e., that $A$ has potentially good reduction at $\mathfrak p$. Since the $n$-torsion is rational over $K$ we conclude, as in Lemma 9.3.1 in \cite{A-HV}, that $A$ has good reduction at $\mathfrak p$.
\end{rem} 
\section{Extension of the Kuga-Satake morphisms in positive characteristic}
The following two sections contain the main results of this note. We show that the Kuga-Satake morphism from Definition \ref{KSLevStrMor} extends in positive characteristic. In this way we give a partial answer to the question posed in the beginning of the chapter.

In Section \ref{AbstractExtension} we prove an abstract extension result concerning morphisms from smooth schemes into $\Av_{g,d',n}$. Then we use this in the next section to show that $ f^{ks}_{d,a,n,\gamma,E_n}$ extends over an open part of $\Spec(\O_{E_n})$. 
\subsection{The extension result}\label{AbstractExtension}
Let us fix the following notations we will use in this section:
\begin{itemize}
\item $R$ will be a discrete valuation ring of mixed characteristic $(0,p)$
where $p > 2$. Denote by $\eta$ and $s$ the generic and the special points of
$\Spec(R)$, respectively. Further, let $K$ be the fraction field of $R$ 
and $k$ will denote the residue field of $R$;
\item $U$ will be a smooth scheme over $R$;
\item We fix three natural numbers $g,d'$ and $n \geq 3$ and denote by $\Av$
the moduli stack $\Av_{g,d',n, R}$ of $g$-dimensional abelian varieties with
polarization of degree ${d'}^2$ and Jacobi level $n$-structure over $R$. We will assume that $p$ does not divide $d'n$.
\item Assume given a morphism $f_\eta \colon U_\eta \lr \Av_\eta$.
\end{itemize}
We are interested in extending the morphism $f_\eta$ over $R$. Of course, in general one cannot expect to be able to do this without further assumptions on $f_\eta$ and $U$. We will list some conditions below which, if satisfied, will guarantee the
existence of an extension of $f_\eta$.  
\begin{assmt}\label{extassumpt}
Let $x$ be a point on the special fiber $U_s$ of $U$, let $\O_{U,x}$
be its local ring and denote by $L$ the  field of fractions of $\O_{U,x}$. Then the morphism
$f_\eta \colon \Spec(L) \lr \Av_\eta$ extends to a morphism $\tilde f \colon
\Spec(\O_{U,x}) \lr \Av$.
\end{assmt}
We will show that if this assumption is fulfilled for certain points of $U$ then the morphism does extend over $R$. More precisely we have
\begin{prp}\label{Ext} 
Let $U$ be a smooth scheme over $R$ and let $f_\eta
\colon U_\eta \lr \Av_\eta$ be a morphism. Assume that the total ramification
index $e$ of $R$ satisfies $e < p - 1$ and that all generic points of the
special fiber $U_s$ of $U$ satisfy Assumption \ref{extassumpt}. Then $f_\eta$ extends uniquely to a morphism $f \colon U \lr
\Av$ over $R$.
\end{prp}
\begin{proof} We will divide the proof into several steps.
\newline
\newline
{\bf Step 1:} We will prove first that if the morphism extends then the extension is unique.
\begin{lem}\label{uniqueext}
Let $V$ be a scheme over $R$ and suppose given a morphism $f_\eta \colon
V_\eta \lr \Av_\eta$.
Assume that it extends to a morphism $f_{V'} \colon V' \lr \Av$ over a dense open subscheme $V'$ of $V$. Then this extension is unique.
\end{lem}
\begin{proof}
This boils down to the fact that $\Av$ is separated over $R$. Assume that there
exist two morphisms, say $F_1$ and $F_2$ extending $f_\eta$ over $V$. Consider the morphism $(F_1,F_2) \colon V \lr \Av \times_R \Av$. The locus where $F_1 = F_2$ is the pull-back $(F_1,F_2)^{-1}\Delta_\Av$ of the diagonal $\Delta_\Av \subset \Av \times_R \Av$ which is closed as $\Av$ is separated over $R$. Hence we conclude that $F_1 = F_2$ on $V'$. 
\end{proof}
In particular, we conclude that if $f$ extends over $U$, then this extension is unique.
\newline
\newline
{\bf Step 2:} There exists a maximal open subscheme $V$ of $U$ such that
$f_\eta$ extends to a morphism $f_V \colon V \lr \Av$. Indeed, if $f$ extends over
two open subschemes $V_1$ and $V_2$ of $U$, then by Lemma \ref{uniqueext} above
those two extensions agree on $V_1 \cap V_2$. Hence we can glue the two morphism
and get a morphism $f_{V_1 \cup V_2} \colon V_1 \cup V_2 \lr \Av$. This shows
that one can take $V$ to be the union of all open subschemes $V'$ of $U$ such
that $f_\eta$ extends to a morphism $f_{V'} \colon V' \lr \Av$.
\newline
\newline
{\bf Step 3:}
Consider the graph $\G_\eta$ of $f_\eta \colon U_\eta \lr \Av_\eta$ in
$U_\eta \times \Av_\eta$. Take the flat extension $\bar \G$ of $\G_\eta$
over $R$ i.e., the closure of $\G_\eta$ in $U \times \Av$. We have the
projection map $pr_1 \colon \bar \G \lr U$. Let $\{U_s^i\}_{\{i \in I\}}$ be the set of connected components of the special fiber $U_s$ and for each $i$ let $U^{(i)}$ be the open subscheme $U \setminus \bigl(\cup_{j\ne i}U^j_s\bigr)$. We look at the set
\begin{displaymath}
U^{(i)}_{\rm ft} := \{ p\in U^{(i)}\ |\ pr_1\ \text{is\ flat\ at\ all\ points\ of}\ pr_1(p) \}.
\end{displaymath}
which is open in $U^{(i)}$ (\S 8, Proposition 8.9.4 in EGA IV (\cite{EGA})). We will call this set, with its induced scheme structure the maximal subscheme of $U^{(i)}$ over which this projection map is flat. Let $V_{\rm ft}$ be the open subscheme $\bigcup_i U^{(i)}_{\rm ft}$ of $U$.
\begin{lem}
The maximal open subscheme $V$ given in Step 2 over which $f_\eta$ extends is equal the open subscheme $V_{\rm ft}$ of $U$ over which the morphism $pr_1 \colon
\bar \G \lr U$ is flat.
\end{lem}
\begin{proof}
Note first that the scheme $V$ from Step 2 is contained in $V_{\rm ft}$. Indeed, the morphism $pr_1\colon pr_1^{-1}(V) \lr V$ is an isomorphism, as $pr_1^{-1}(V) \subset U \times \Av$ is the graph of $f_V$, and therefore it is flat.

Since $pr_1 \colon pr^{-1}_1(V_{\rm ft}) \lr V_{\rm ft}$ is flat the dimension
of the fibers is constant. It is zero on the generic fiber hence this morphism
is quasi-finite. Using Zariski's Main Theorem (\S 8, Theorem
8.12.6 in EGA IV (\cite{EGA})) one can factor $pr_1|_{V_{\rm ft}}$ as an open immersion and a
finite morphism $pr^{-1}_1(V_{\rm ft}) \lr \tilde V \lr V_{\rm ft}$. But
generically, over every connected component of $V_{\rm ft}$, the degree of the finite morphism is one, hence it is one everywhere. This means that $pr_1 \colon pr^{-1}_1(V_{\rm ft}) \lr V_{\rm ft}$ is an
isomorphism. Hence one can extend $f_\eta$ on $V_{\rm ft}$ using the second
projection map $pr_2 \colon \bar \G \lr \Av$. Therefore we get the other
inclusion $V_{\rm ft} \subset V$.
\end{proof}
{\bf Step 4:} We will show that the open subscheme $V$ given is Step 2
contains all generic points of $U_s$. We need some auxiliary results.
\begin{lem}\label{Lemm3.5}
If $x \in U$ satisfies Assumption \ref{extassumpt}, then $x \in V$.
\end{lem}
\begin{proof}
According to Step 3 we have to show that $pr_1 \colon \bar \G \lr U$ is flat at
$x$.
\begin{claim}\label{claim1}
Let $X$ be a scheme over $R$ and let $\G_\eta \subset X_\eta$
be a closed subscheme. Take the flat extension $\bar \G$ of $\G_\eta$. Let $i
\colon Y \lr X$ be a flat morphism over $R$ and set $\Delta_\eta \colon =
i^*\G_\eta$. Let $\bar \Delta$ be the flat extension of $\Delta_\eta$ in $Y$.
Then one has that $\bar \Delta = i^*\bar \G$.
\end{claim}
\begin{proof}
Since flatness is stable under base change we have that $i^*\bar \G \lr \bar \G$ is
flat. Hence $i^*\bar \G \lr R$ being the composite $i^*\bar \G
\lr \bar \G \lr R$ of two flat maps is also flat. Moreover, by the very definitions we have
that $(i^*\bar \G)_\eta = i^*\G_\eta = \Delta_\eta$. Therefore by uniqueness
of the flat extension (see \S2, Proposition 2.8.5 in EGA IV (\cite{EGA})) we conclude that
$\bar \Delta = i^*\bar \G$.
\end{proof}
For a scheme $X$ denote by $|X|$ its underlying topological space.
\begin{claim}\label{claim2}
Let $X$ and $T$ be two schemes and $h \colon T \lr X$ be a
morphism. Take a point $x \in |X|$, let $i \colon \Spec(\O_{X,x}) \lr X$ be the morphism associated to $x$ and let $i^*h \colon i^*T \lr \Spec(\O_{X,x})$ be the pull-back map. If $i^*h$ is flat, then $h$ is flat above $x$ i.e., it is flat at all points $t \in h^{-1}(x)$.
\end{claim}
\begin{proof}
If $t \in h^{-1}(x)$, then $\O_{T,t} \cong \O_{i^*T,t}$ as $\O_{X,x}$-algebras and hence we obtain the result in the claim.
\end{proof}
Let us go back to the proof of the Lemma \ref{Lemm3.5}. We have the following diagram:
$$
\xymatrix{
i^*\bar \G \ar[rr] \ar[d] & & \G \ar[d] \\
\Spec(\O_{U,x}) \times \Av \ar[d] \ar[rr]^{\rm flat} & & U \times \Av \ar[d] \\
\Spec(\O_{U,x}) \ar[dr]_{\rm flat} \ar[rr]^{\rm flat} & & U \ar[dl]^{\rm smooth}
\\
 & \Spec(R). &
}
$$
Let $\a = \Spec(L)$ where $L$ is the field of fractions of $\O_{U,x}$. Consider
the point $\Delta_\eta = (\a, f_\eta(\a))$ on $U_\eta \times \Av_\eta$. Then by
Claim \ref{claim1} applied to 
\begin{align*} 
 X &= \Spec(\O_{U,x}) \times \Av \\
 Y &= U \times \Av\ \text{over}\ R \\
\G_\eta  &= {\rm the\ graph\ of}\ f_\eta
\end{align*}
we conclude that $i^*\bar\G = \bar \Delta$, where $\bar\Delta$ is the flat closure over $R$.

Since Assumption \ref{extassumpt} holds for the point $x$ the map $f_\eta \colon \Spec(L) \lr \Av_\eta$ extends to a morphism $\tilde f
\colon \Spec(\O_{U,x}) \lr \Av$ and we see that $\bar\Delta$ is the graph of $\tilde
f$. In particular we have that $\bar \Delta \cong \Spec(\O_{U,x})$ hence it is flat over
$\Spec(\O_{U,x})$. If we apply Claim \ref{claim2} with $T = \bar\G$ and $X = U$ we get
that $pr_1 \colon \bar \G \lr U$ is flat at $x$. Therefore by Step 3 we conclude
that $x \in V$. This finishes the proof of the Lemma.
\end{proof}
Since all generic points of the special fiber $U_s$ satisfy Assumption \ref{extassumpt} we conclude that they are contained in $V$.
\newline
\newline 
{\bf Step 5:} By Step 4 there exists an open dense subscheme $V$ of $U$,
containing the generic points of the connected components of the special fiber
$U_s$ over which the morphism $f_\eta \colon U_\eta \lr \Av_\eta$ extends to a
morphism $f_V \colon V \lr \Av$ over $R$. It corresponds to a polarized abelian
scheme with level $n$-structure $(A_V, \l_V, \a_V)$ over $V$.

As $V$ contains strictly the generic fiber $U_\eta$ and all generic points of
the special fiber $U_s$ we have that $\codim_U U \setminus V \geq 2$. Since $U$
is smooth over $R$ and by assumption $e < p-1$ then by a result of Faltings (Lemma 3.6 in \cite{BM-SV}) one concludes that $A_V \lr V$ extends to an abelian scheme $A \lr U$.

Further, by Lemma \ref{polext} the polarization $\l_V$ extends to a polarization
$\l \colon A \lr A^t$. Since $p$ does not divide $n$, the level $n$-structure $\a_V$ extends uniquely to a level $n$-structure $\a$ on $(A,\l)$. Hence we get a polarized abelian scheme $(A, \l, \a)$ extending $(A_V, \l_V, \a_V)$. This corresponds to a morphism $f \colon U \lr \Av$ extending $f_\eta$. 
\end{proof}
\begin{rem}\label{ProblemExtKSQ}
We will apply Proposition \ref{Ext} to show that the Kuga-Satake morphism constructed in Proposition \ref{KSCMapN} extends over an open part of $\Spec(\O_{E_n})$, where $\O_{E_n}$ is the ring of integers in $E_n$. One might want to use the same line of thoughts and try to extend the Kuga-Satake morphism $f_{d,a,\Kg^{\rm a}, \Q}^{ks}$ defined in Proposition \ref{KSCmap} over an open part of $\Spec(\Z)$. The problem which one comes up with is to carry on Step 3 in this situation. One can define an equivalence of the closure $\bar \Gamma$ of $\Gamma$. In general, the morphism $\bar \Gamma_{\rm ft} := pr^{-1}(V_{\rm ft}) \lr V_{\rm ft}$ might not be representable so one cannot use Zariski's Main Theorem (Theorem 16.5 in \cite{L-MB}).
\end{rem}
\subsection{Extension of the Kuga-Satake morphisms}\label{KSPosChar}
In this section we will use the notations established in \S \ref{KSMorphismC}. In particular, we fix two natural numbers $d$ and $n$ and let us suppose that $n \geq 3$. Let $\gamma$ be a set-theoretic section of the homomorphism $\Gal(E_n/\Q) \lr \Gal(E_X/\Q)$. We will show below that the Kuga-Satake morphism $f^{ks}_{d,a,n,\gamma,E_n}$ extends over an open part of $\Spec(\O_E)$ where $\O_{E_n}$ the ring of integers in $E_n$. 

\begin{thm}\label{abstarext}
Let $d,n \in \N$, $n\geq 3$ and suppose that $a\in C^+$ satisfies Assumption \ref{PolAssumpt}. Then the Kuga-Satake morphism $f^{ks}_{d,n,a,\gamma,E_n} \colon \Fk_{2d,n^{\rm sp},E_n} \lr \Av_{g,d',n,E_n}$ extends uniquely to a morphism
\begin{displaymath}
 f^{ks}_{d,a,n,\gamma} \colon \Fk_{2d,n^{\rm sp}, \O_{E_n}[1/N]} \lr \Av_{g,d',n,\O_{E_n}[1/N]}
\end{displaymath}
where $N = 2dd'nl$ and $l$ is the product of the prime numbers $p$ whose ramification index $e_p$ in $E_n$ is $\geq p-1$.
\end{thm}
\begin{proof}
Let us first shorten the notations a bit by setting $\Fk$ to be $\Fk_{2d,n^{\rm sp}, \O_{E_n}[1/N]}$ and $\Av$ to be $\Av_{g,d',n,\O_{E_n}[1/N]}$. Let $\pi \colon U \lr \Fk$ be an atlas of $\Fk$ (i.e., an \'etale surjective morphism) over $\O_{E_n}[1/N]$. We may assume that the pull-back of the universal family of polarized K3 surfaces to $U$ is a K3 scheme. The map $f^{ks}_{d,a,n,\gamma,E_n} \colon \Fk_{E_n} \lr \Av_{E_n}$ defines a morphism $f_{E_n} = f^{ks}_{d,a,n,\gamma, E_n} \circ \pi_{E_n} \colon U_{E_n} \lr \Av_{E_n}$. We will fist extend $f_{E_n}$ to a morphism over $\O_{E_n}[1/N]$ and then using a descent argument show that it comes from a morphism $f^{ks}_{d,a,n,\gamma} \colon \Fk_{d,n^{\rm sp}, \O_{E_n}[1/N]} \lr \Av_{g,d',n,\O_{E_n}[1/N]}$.

Let $\mathfrak p$ be a prime ideal of $E_n$ not dividing $N$ and let $R = \O_{E_n, (\mathfrak p)}$ be the localization of $\O_{E_n}$ at $\mathfrak p$. As before $\{s,\eta\}$ will be the special and the generic points of $\Spec(R)$. In order to apply Proposition \ref{Ext} to $f_{E_n}$ and $U_R$, which we will denote by $U^{\mathfrak p}$, we have to show that all generic points of the special fiber $U^{\mathfrak p}_s$ of $U^{\mathfrak p}$ satisfy Assumption \ref{extassumpt}.

Let $x \in |U^{\mathfrak p}_s|$ be a generic point. Then $\O_{U^{\mathfrak p},x}$ is a discrete valuation ring with a maximal ideal $\mathfrak{m}_x$. Let us denote its field of fractions by $L$. Taking the pull-back of the universal family of polarized K3 surfaces with spin level $n$-structures via the canonical morphism $\Spec(\O_{U^{\mathfrak p}, x}) \lr U^{\mathfrak p}$ we obtain a K3 scheme $(X \lr \Spec(\O_{U^{\mathfrak p}, x}), \l, \nu)$. Then $f_{E_n}$ gives a morphism $\Spec(L) \lr \Av_{E_n}$ and let $(A, \mu,\e)$ be the corresponding abelian variety over $L$. It is the Kuga-Satake abelian variety associated to the generic fiber of $(X \lr \Spec(\O_{U^{\mathfrak p}, x}), \l, \nu)$. We can apply Remark \ref{GRedProperty} (or alternatively by Lemma 9.3.1 in \cite{A-HV} we conclude that $A$ has potentially good reduction and as the $n$-torsion points are $L$-rational, then $A$ has good reduction) to see that the abelian variety $A$ has good reduction at $\mathfrak{m}_x$. In other words the N\'eron model of $A$ over $\O_{U^{\mathfrak p},x}$ is an abelian scheme. By Lemma \ref{pol}, the polarization $\l$ extends uniquely over $\O_{U^{\mathfrak p},x}$ and as $\mathfrak p$ does not divide $n$, the level $n$-structure extends uniquely, as well. Hence the morphism $\Spec(L) \lr \Av_{E_n}$ extends to a morphism $\Spec(\O_{U^{\mathfrak p}, x}) \lr \Av_R$. Therefore, by Proposition \ref{Ext} applied to $U^{\mathfrak p}$, $f_{E_n}$ and $R=\O_{E_n, (\mathfrak p)}$ one can extend $f_{E_n} \colon U_{E_n} \lr \Av_{E_n}$ to a morphism $f_{\mathfrak p} \colon U^{\mathfrak p} \lr \Av_R$.

The morphism $f_{E_n}$ can be extended uniquely over $\O_{E_n,(\mathfrak p)}$ for any $p$ not dividing $N$. Hence we conclude that it extends uniquely to a morphism $f \colon U \lr \Av$ over $\Z[1/N]$.

We are left to show that $f$ descends to $\Fk_{2d,n^{\rm sp},\O_{E_n}[1/N]}$. By Proposition 1.4 in \S 1, Ch. II of \cite{Knu-AS} one has the following exact sequence
$$
\xymatrix{
 0 \ar[r] & \Hm_S(\Fk,\Av ) \ar[r]^{\a^*} & \Hm_S(U, \Av)
 \ar@<2pt>[r]^{pr^*_1} \ar@<-2pt>[r]_{pr^*_2} & \Hm_S(U',\Av)
}
$$
where $U' = U \times_{\Fk} U$. Note that both $pr^*_1(f)$ and $pr^*_2(f)$ are extensions of the morphism $pr^*_1 \circ \pi^* (f^{ks}_{d,a,n,\gamma, E_n}) = pr^*_2 \circ \pi^*(f^{ks}_{d,a,n,\gamma,E_n})$ over $\O_{E_n}[1/N]$. Since $U'$ is a smooth scheme over $\O_{E_n}[1/N]$ (Definition 1.1 in \cite{Knu-AS}) and $\Av$ is separated just like in Lemma \ref{uniqueext} we conclude that such an extension is unique. Hence we one has that $pr^*_1(f) = pr^*_2(f)$ and therefore by the above exact sequence $f$ comes from a morphism
$$
 f^{ks}_{d,a,n,\gamma} \colon \Fk_{2d,n^{\rm sp}, \O_{E_n}[1/N]} \lr \Av_{g,d',n,\O_{E_n}[1/N]}
$$
over $\O_{E_n}[1/N]$.
\end{proof}
We end this section with a few remarks concerning the Kuga-Satake morphism in mixed characteristic.
\begin{rem}
In the proof of Proposition \ref{Ext} we used a result of Faltings to show that certain morphisms extend in positive characteristic. This is really an essential step of our strategy for defining Kuga-Satake abelian varieties in positive characteristic. In Theorem \ref{abstarext} we have to exclude the primes $p$ for which the ramification index $e_p$ is $\geq p-1$, as Lemma 3.6 in~\cite{BM-SV} does not hold for these primes. See Section 6 in \cite{dJO} and Section 3.4 in~\cite{BM-SV}.
\end{rem}
\begin{rem}
We use the notations of Theorem \ref{abstarext}. Suppose $k$ is a field of characteristic $p$ such that $p$ does not divide $N$ and let $R$ be a discrete valuation ring of mixed characteristic $(0,p)$ with field of fractions $k$. Then to every primitively polarized K3 surface with a spin level $n$-structure $(X,\l,\nu)$ over $k$ we associate via $f^{ks}_{d,a,n,\gamma}$ a polarized abelian variety with level $n$-structure $(A,\mu,\e)$ over $k$. We will call $A$ the \emph{Kuga-Satake abelian variety} associated to $(X,\l,\nu)$. Further, if $(X_1,\l_1,\nu_1)$ and $(X_2,\l_2,\nu_2)$ are two lifts of $(X,\l,\nu)$ over $R$, then the special fibers of $(A_i,\mu_i,\e_i) := f^{ks}_{d,a,n,\gamma}\bigl((X_i,\l_i,\nu_i)\bigr)$, for $i =1,2$ are the same.
\end{rem}
\begin{rem} 
In characteristic zero one can show that the image $f^{ks}_{d,n,a,\gamma,E_n}(\Fk_{2d,n^{\rm sp},E_n})$ in $\Av_{g,d',n,E_n}$ is locally closed. Indeed, as we saw in Proposition 2.5 the period map is open and the morphisms $\beta_{(\Kg^{\rm sp}_n,\Lambda_n)}$ and $\delta_\gamma$ involved in the construction of $f^{ks}_{d,n,a,\gamma,E_n}$ are finite (see Definition \ref{KSLevStrMor}). It is interesting to know if the same holds in mixed characteristic. This question is directly connected to the existence of an analogue of the N\'eron-Ogg-Shafarevich criterion for potentially good reduction of K3 surfaces. To our knowledge, in general, this is still an open problem.
\end{rem}
\section{Applications}
We end this note with some applications of the existence of Kuga-Satake morphisms in mixed characteristic. In Section \ref{CohGrCompSect} we show that the \'etale cohomology relations \S 6.6.1 \cite{D-K3} and in Definition 4.5.1 in \cite{A-HV} hold for the Kuga-Satake abelian varieties defined in \S \ref{KSPosChar}. Then, in Section \ref{CanLiftSection}, we study the behavior of $ f^{ks}_{d,a,n,\gamma}$ at ordinary points. Suppose that $k$ is a finite field of characteristic $p$ where $p$ does not divide $N$ (cf. Theorem \ref{abstarext}) and let $(X,\L,\nu) \in \Fk_{2d,n^{\rm sp},\F_p}(k)$ be an ordinary point. We will prove that the canonical lift $(X^{\rm can},\L,\nu)$ over $W(k)$ is mapped to the canonical lift $(A^{\rm can},\mu,\e)$ of $(A,\mu,\e) =  f^{ks}_{d,a,n,\gamma}\bigl((X,\L,\nu)\bigr)$.

\subsection{Cohomology groups}\label{CohGrCompSect}
Let $d$ and $n$ be two natural numbers and suppose further that $n \geq 3$. With the notations as in Sections \ref{KSMorphismC} and \ref{KSPosChar} let $a \in C^+$ be an element satisfying Assumption \ref{PolAssumpt} and let $\gamma$ be a set-theoretic section of the homomorphism $\Gal(E_n/\Q) \lr \Gal(E_X,\Q)$. Then we have a Kuga-Satake morphism
$$
 f^{ks}_{d,a,n,\gamma} \colon \Fk_{2d,n^{\rm sp}, \O_{E_n}[1/N]} \lr \Av_{g,d',n,\O_{E_n}[1/N]}
$$
where $N = 2dd'nl$ and $l$ is the product of the prime numbers $p$ whose ramification index $e_p$ in $E_n$ is $\geq p-1$.

Let $k$ be a field of characteristic $p$ and suppose given a $k$-valued point $(X,\l,\nu) \in \Fk_{2d,n^{\rm sp}, \O_{E_n}[1/N]}$. Denote by $(A,\mu,\beta_n)$ the polarized Kuga-Satake abelian variety with level $n$-structure $f^{ks}_{d,a,n,\gamma}((X,\l,\nu))$.
\begin{lem}
With the notations as above one has and isomorphism of $\Gal(\bar k/k)$-modules
$$
C^+\bigl(P^2_{\rm et}(X_{\bar k}, \Z_l(1))\bigr) \cong \End_{C^+}(H^1_{\rm et}(A_{\bar k},\Z_l))
$$
for any $l \ne p$.
\end{lem}
\begin{proof}
Let $(\mathcal X, \l, \nu)$ be a lift of $(X,\l,\nu)$ over $W(k)$ (which exists because $\Fk_{2d,n^{\rm sp}, \O_{E_n}[1/N]}$ is smooth over $\O_{E_n}[1/N]$) and let $(\Av, \mu, \e)$ be the Kuga-Satake variety $f^{ks}_{d,a,n,\gamma}\bigl((\mathcal X, \l, \nu)\bigr)$. By Corollary \ref{ChoCompRes} we have an isomorphism of $\Gal(\bar K /K)$-modules
$$
C^+\bigl(P^2_{\rm et}(\X_{\bar K}, \Z_l(1))\bigr) \cong \End_{C^+}(H^1_{\rm et}(\Av_{\bar K},\Z_l))
$$
for any $l$. Hence if $l\ne p$ one can apply the smooth base change theorem for \'etale cohomology to prove the claimed isomorphism.
\end{proof}
\begin{rem}
Note that one can use this isomorphism in case $k = \F_q$ to compute the Newton polygon of $A$ in terms of the Newton polygon of $X$. For instance one can see that if $X$ is ordinary then $A$ is also ordinary. We refer to Proposition 2.5 in \cite{Nyg-TC} for a proof.
\end{rem} 

\subsection{Canonical lifts of ordinary K3 surfaces}\label{CanLiftSection}
Let $k$ be a perfect field of characteristic $p>0$ and let $W(k)$ be the ring of Witt vectors. Suppose given an ordinary K3 surface $X_0$ over $k$. Denote by $\mathcal X / S$ the universal deformation of $X_0$ over $W(k)$. We know that $S$ is formally smooth of dimension 20 (see Corollary 1.2 in \cite{Del-K3}).

In Section IV of \cite{A-M}, Artin and Mazur define the \emph{enlarged formal Brauer group} $\Psi_{X_0}$ of $X_0$ which is a $p$-divisible group over $k$ and such that its connected component $\Psi^0_{X_0}$ is $\hat Br(X_0)$. With the notations of Section 4.1 in \cite{Riz-MK3}, let $R \in \underline A$ be a local artinian ring with residue field $k$ and let $(X \lr \Spec(R),\phi)$ be a deformation of $X_0$ over $R$. Then the enlarged formal Brauer group $\Psi_X$ over $\Spec(R)$ exists. It is a $p$-divisible group over $\Spec(R)$ and the isomorphism $\phi$ induces an isomorphism of $p$-divisible groups over $\Spec(k)$
\begin{displaymath}
 \phi_{Br} \colon \Psi_X \otimes_R k \lr \Psi_{X_0}.
\end{displaymath}   
In other words $\Psi_X$ is a lifting of $\Psi_{X_0}$ over $R$. Let
\begin{displaymath}
\DefF_{\Sch}(X_0) \colon \underline A \lr {\rm Sets}
\end{displaymath}
be the covariant deformation functor defined in Section 4.1 in \cite{Riz-MK3} and let 
\begin{displaymath}
{\rm DefBr}_{X_0} \colon \underline A \lr {\rm Sets}
\end{displaymath}
be the covariant functor
\begin{displaymath} 
\begin{matrix}
{\rm DefBr}_{X_0}(R) = & \bigl\{{\rm isomorphism\ classes\ of\ pairs}\  (G,\phi)\ {\rm where}\ G\ {\rm is\ a}\\ 
& p-{\rm divisible\ group\ over}\ R\ {\rm and}\ \phi \colon G \otimes_R k \cong \Psi_{X_0}\bigr\}.
\end{matrix} 
\end{displaymath}
We have the following Serre-Tate theory for ordinary K3 surfaces:
\begin{thm}[Nygaard]\label{STK3}
For any $R \in \underline A$ the map
\begin{displaymath}
\DefF_{\Sch}(X_0)(R) \lr {\rm DefBr}_{X_0}(R)
\end{displaymath}
defined by
\begin{displaymath}
 (X\lr \Spec(R),\phi) \mapsto (\Psi_X, \phi_{Br})
\end{displaymath}
is a bijection.
\end{thm}
\begin{proof}
For a proof we refer to Theorem 1.1 in \cite{Nyg-TC}.
\end{proof}
Let $G$ be a lifting of $\Psi_{X_0}$ over $R$. Since height one groups are rigid, we have precisely one lifting $G^0_R$ of $\hat Br(X_0) = \Psi^0_{X_0}$ to $\Spec(R)$. Similarly, \'etale groups are also rigid so there is a unique lift $G^{\rm et}_R$ of $\Psi^{\rm et}_{X_0}$ to $\Spec(R)$. So we for any lifting $G$ of $\Psi_{X_0}$ to $\Spec(R)$ we have en exact sequence
\begin{displaymath}
 0 \lr G^0_R \lr G \lr G^{\rm et}_R \lr 0
\end{displaymath}
lifting
\begin{displaymath}
 0 \lr \hat Br(X_0) \lr \Psi_{X_0} \lr \Psi^{\rm et}_{X_0} \lr 0
\end{displaymath}
over $\Spec(R)$. 

If we consider the trivial extension $G = G^0_R\times G^{\rm et}_R$, then by Theorem \ref{STK3} above there is a unique lifting $X_R^{\rm can}$ of $X_0$ over $\Spec(R)$ such that $\Psi_{X_R^{\rm can}} = G^0_R \times G^{\rm et}_R$. For any $n \in N$ taking $R = W_n$ we obtain a lifting $X_n = X_{W_n}^{\rm can}$. The projective system $\{X_n\}$ defines a proper flat formal scheme $\{X_n\}$ over $\Spf(W)$. It is algebrizable and defines a K3 scheme $X^{\rm can}$ over $\Spec(W)$ called the \emph{canonical lift} of $X_0$. Every line bundle of $X_0$ lifts uniquely to a line bundle on $X^{\rm can}$. For a proof of these facts we refer to Proposition 1.8 in \cite{Nyg-TC}.

With the notations established in Section \ref{KSMorphismC} let $d$ and $n \geq 3$ be two natural numbers, and let $a \in C^+$ be an element satisfying Assumption \ref{PolAssumpt}. Choose a set-theoretic section $\gamma$ of the homomorphism $\Gal(E_n/\Q) \lr \Gal(E_X/\Q)$ so that we have a Kuga-Satake morphism
\begin{displaymath}
 f^{ks}_{d,a,n,\gamma} \colon \Fk_{2d,n^{\rm sp}, \O_{E_n}[1/N]} \lr \Av_{g,d',n,\O_{E_n}[1/N]}
\end{displaymath}
where $N = 2dd'nl$ and $l$ is the product of the prime numbers $p$ whose ramification index $e_p$ in $E_n$ is $\geq p-1$. Let $k = \F_q$ be a finite field and suppose given an ordinary point $(X_0,\L_0,\nu_0) \in \Fk_{2d,n^{\rm sp},\O_{E_n}[1/N]}(k)$ (in particular $p$ does not divide $N$). Denote by $(X^{\rm can},\L)$ the canonical lift of $X_0$ to $W$. The spin level $n$-structure $\nu_0$ also lifts uniquely to a spin level $n$-structure on $X^{\rm can}$ as $p$ does not divide $N$. Denote by $(A^{ks},\mu,\e)$ the abelian scheme $f^{ks}_{d,a,n,\gamma}\bigl((X^{\rm can},\L,\nu)\bigr)$ over $\Spec(W)$ and let $(A_0,\mu_0,\e_0)$ be the triple $(A^{ks},\mu,\e) \otimes k$ over $k$. 

The following result was suggested to us by B. Moonen.
\begin{prp}\label{CanLift}
The abelian scheme $A^{ks}$ is the canonical lift of $A_0$ over $\Spec(W)$.
\end{prp}
\begin{proof}
By Theorem 2.7 in \cite{Nyg-TC} we know that, after a base change $R' \lr R$, the abelian scheme $A^{ks}$ is isogenous to the canonical lift $A^{\rm can}$ of $A_0$. Hence we conclude that $A^{ks}$ is a quasi-canonical lift.

Let 
\begin{displaymath}
 {\rm Def}_{(A_0,\mu_0)} \colon \underline A \lr {\rm Sets}
\end{displaymath}
be the covariant functor
\begin{displaymath}
\begin{matrix}
{\rm Def}_{(A_0,\mu_0)}(R) = &  \{{\rm isom.\ classes\ of\ polarized\ abelian\ schemes\ (A,\mu, \phi)\ over}\ $R$ \\
  &  {\rm and\ an\ isomorphism}\ \phi \colon (A,\mu)\otimes_R k \lr (A_0,\mu_0) \}.
\end{matrix}
\end{displaymath}
This functor is representable by a formal smooth scheme $\mathfrak A_{(A_0,\mu_0)}$ which has a structure of a formal torus. For details we refer to Theorems 1.2.1 and 2.1 in \cite{Katz-LMod} and Ch. III, \S 1 of \cite{BM-Thesis}.

The lift $A^{ks}/W$ defines a point $s \in \mathfrak A_{(A_0,\mu_0)}(W)$. Just like in Lemma 1.5 in Chapter~{III}, \S~1 in~\cite{BM-Thesis} we conclude that since $A^{ks}$ is a quasi-canonical lift then $s$ is a torsion point. As $\mathfrak A_{(A_0,\mu_0)}(W)$ is $l$-divisible for all $l \ne p$ we have that $s^{p^m} = 1 \in \mathfrak A_{(A_0,\mu_0)}(W)$. But $s$ is defined over $W$ which is unramified and any $p^m$ torsion point is defined over ramified rings unless $m = 0$. Hence, $s = 1$ which corresponds to the canonical lift in the \emph{Serre-Tate coordinates}. Therefore, $A^{ks}$ is the canonical lift of its special fiber.
\end{proof}

Let $\mathfrak p$ be a prime ideal of $\O_{E_n}$ which does not divide $N$ and let $k := \O_{E_n}/\mathfrak p$ be its residue field. It is a finite field and let $p$ be its characteristic. Let $R$ be the localization of $\O_{E_n}$ at $\mathfrak p$. It is a discrete valuation ring with a residue field $k$. Let $\Fk^{(2)}_{2d,n^{\rm sp},k}$ be the non-ordinary locus of $\Fk_{2d,n^{\rm sp},k}$ (see Section 1 in \cite{JR-HSK3}). It is a closed subspace of $\Fk_{2d,n^{\rm sp},k}$ and we consider
\begin{displaymath}
\Fk^{\rm ord}_{2d,n^{\rm sp},R} := \Fk_{2d,n^{\rm sp},R} \setminus \Fk^{(2)}_{2d,n^{\rm sp},k}
\end{displaymath}
which is an open subspace of $\Fk_{2d,n^{\rm sp},R}$.
\begin{cor}
The restriction of the Kuga-Satake morphism
\begin{displaymath}
 f^{ks}_{d,a,n,\gamma,R} \colon \Fk^{\rm ord}_{2d,n^{\rm sp},R} \lr \Av_{g,d',n,R}
\end{displaymath}
is quasi-finite.
\end{cor}
\begin{proof}
Note first that $f^{ks}_{d,a,n,\gamma,E_n}$ is a quasi-finite morphism. Indeed, by construction (see Definition \ref{KSLevStrMor}) we have that $f^{ks}_{d,a,n,\gamma,E_n} = \beta_{(\Kg^{\rm sp}_n,\Lambda_n)} \circ \delta_\gamma \circ j_{d,\Kg^{\rm a},E_n}$ where 
\begin{displaymath}
 j_{d,\Kg^{\rm a}_n,E_n} \colon \Fk_{2d,n^{\rm sp},E_n} \lr Sh_{\Kg^{\rm sp}_n}(G,\Omega^\pm)_{E_n}
\end{displaymath}  
is an \'etale morphism of noetherian schemes and 
\begin{displaymath}
 \beta_{(\Kg^{\rm sp}_n,\Lambda_n)} \circ \delta_\gamma \colon Sh_{\Kg^{\rm sp}_n}(G,\Omega^\pm)_{E_n} \lr \Av_{g,d',n,E_n}
\end{displaymath}  
is a quasi-finite morphism. Therefore $f^{ks}_{d,a,n,\gamma,E_n}$ is a quasi-finite morphism. To finish the proof we have to show that for any $\bar k$-valued point $y \in \Av_{g,d',n,k}(\bar k)$ there are only finitely may $\bar k$-valued points $x \in \Fk^{\rm ord}_{2d,n^{\rm sp},k}(\bar k)$ such that $f^{ks}_{d,a,n,\gamma,k}(x) = y$.

Suppose that $(X_1,\l_1,\nu_1)$ and $(X_2,\l_2,\nu_2)$ are two ordinary K3 surfaces over a finite field $L \subset k$ in $\Fk_{2d,n^{\rm sp},k}(L)$ such that 
\begin{displaymath}
 f^{ks}_{d,a,n,\gamma,k}\bigl((X_1,\l_1,\nu_1)\bigr) =  f^{ks}_{d,a,n,\gamma,k}\bigl((X_2,\l_2,\nu_2)\bigr) = (A,\mu,\e).
\end{displaymath}
Taking a finite extension of $L$, if needed, we may assume that $\l_1$ and $\l_2$ are classes of ample line bundles $\L_1$ and $\L_2$ on $X_1$ and $X_2$, respectively. Let $(X_1^{\rm can},\L_1,\nu_1)$ and $(X_2^{\rm can},\L_2,\nu_2)$ be the two canonical lifts over $W(L)$. Denote the field of fractions of $W(L)$ by $K$. We have that $(X_1,\L_1,\nu_1) \cong (X_2,\L_2,\nu_2)$ if and only if $(X_1^{\rm can},\L_1,\nu_1)\otimes K \cong (X_2^{\rm can},\L_2,\nu_2)\otimes K$. By Proposition \ref{CanLift} we have that 
\begin{displaymath}
f^{ks}_{d,a,n,\gamma,E_n}\bigl((X_1^{\rm can},\L_1,\nu_1)\otimes K \bigr) =  f^{ks}_{d,a,n,\gamma,k}\bigl((X_2^{\rm can},\L_2,\nu_2)\otimes K\bigr) = (A^{\rm can},\mu,\e)\otimes K
\end{displaymath}
hence we conclude the $f^{ks}_{d,a,n,\gamma,k}$ is quasi-finite from the fact that $f^{ks}_{d,a,n,\gamma,E_n}$ is quasi-finite.
\end{proof}
Combining Theorem 16.5 in \cite{L-MB} and Corollary 6.16 in Ch. II, \S 6 of \cite{Knu-AS} with the preceding corollary we obtain the following result.
\begin{cor}
There exists a scheme $Z$ over $R$, a finite morphism $\pi \colon Z \lr \Av_{g,d',n,R}$ and an open immersion $i \colon \Fk^{\rm ord}_{2d,n^{\rm sp},R} \hookrightarrow Z$ such that  $f^{ks}_{d,a,n,\gamma,R} = \pi \circ i$. Therefore the algebraic space $\Fk^{\rm ord}_{2d,n^{\rm sp},R}$ is a scheme.
\end{cor}
\bibliographystyle{hplain}

\end{document}